\documentclass[leqno,12pt,twoside]{amsart}
\numberwithin{equation}{section}

\usepackage{amsopn}
\usepackage[latin1]{inputenc}
\usepackage[T1]{fontenc}
\usepackage{amsmath, amssymb}
\usepackage[dvips]{graphicx}
\usepackage{epsfig}
\usepackage{psfig}
\usepackage{fullpage}
\usepackage{euscript}
\input xy
\xyoption{all}

\theoremstyle{plain}

	\newtheorem{theorem}{Theorem}[section]
	\newtheorem{theo}{Theorem}

	\newtheorem{Edefinition}[theorem]{Definition}

	\newtheorem{DefProp}[theorem]{Definition -- Proposition}

	\newtheorem{sublemma}[theorem]{Sub-lemma}
	\newtheorem{lemma}[theorem]{Lemma}
	
	\newtheorem{proposition}[theorem]{Proposition}

	\newtheorem{corollary}[theorem]{Corollary}

	\newtheorem{rema}[theorem]{\it Remark}
	\newenvironment{remark}{\begin{rema} \normalfont}{\end{rema}}

	\renewenvironment{proof}
	{\vspace{-0.3cm}\trivlist\item[\hskip \labelsep {\it Proof.
	    }\enskip]}%
	{\unskip\nobreak\hskip 2em plus 1fil\nobreak%
	  $\blacksquare$%
	  \parfillskip=0pt \endtrivlist}
	
	\newenvironment{SketchProof}
	{\vspace{-0.3cm}\trivlist\item[\hskip \labelsep {\it Sketch of proof.
	    }\enskip]}%
	{\unskip\nobreak\hskip 2em plus 1fil\nobreak%
	  $\blacksquare$%
	  \parfillskip=0pt \endtrivlist}

\newcommand{\scal}{\mathrm{scal}}
\newcommand{\vol}{\mathrm{vol}}

\newcommand{\Ricci}{\mathrm{Ricci}}
\newcommand{\Karea}{\operatorname{K-area}}

\newcommand{\RR}{\mathbb{R}}
\newcommand{\ZZ}{\mathbb{Z}}
\newcommand{\CC}{\mathbb{C}}

\newcommand{\NN}{\mathbb{N}}
\newcommand{\QQ}{\mathbb{Q}}

\renewcommand{\d}{\mathrm{d}}
\newcommand{\pd}{\partial}
\newcommand{\dd}[1]{\frac{\d}{\d #1}}
\newcommand{\pdd}[1]{\frac{\pd}{\pd #1}}
\renewcommand{\div}{\mathrm{div}}

\newcommand{\sflow}{\mathrm{sf}}
\renewcommand{\L}{\mathrm{L}}
\newcommand{\D}{\mathrm{D}}
\newcommand{\Dirac}{/ \!\!\!\! D}
\newcommand{\Spin}{/ \!\!\! S}
\newcommand{\tD}{\widetilde{\mathrm{D}}}
\newcommand{\tDirac}{\widetilde{/ \!\!\!\! D}}

\newcommand{\tV}{\widetilde{V}}
\newcommand{\tg}{\widetilde{\mathrm{g}}}
\newcommand{\ind}{\operatorname{ind}}
\newcommand{\tnabla}{\widetilde{\nabla}}
\newcommand{\CB}{\operatorname{CB}}
\newcommand{\spec}{\operatorname{spec}}
\newcommand{\R}{\EuScript{R}}

\renewcommand{\H}{\operatorname{H}}
\newcommand{\ch}{\operatorname{ch}}
\newcommand{\Ah}{\operatorname{\widehat{A}}}

\renewcommand{\deg}{\operatorname{deg}}

\newcommand{\Gl}{\operatorname{Gl}}
\newcommand{\id}{\operatorname{id}}
\newcommand{\rank}{\operatorname{rank}}

\newcommand{\G}{\Gamma}
\newcommand{\Id}{\operatorname{Id}}
\newcommand{\dom}{\operatorname{dom}}
\renewcommand{\Im}{\operatorname{Im}}
\newcommand{\Ker}{\operatorname{Ker}}
\newcommand{\grad}{\operatorname{grad}}
\newcommand{\Tr}{\operatorname{Tr}} 

\renewcommand{\tilde}{\widetilde}
\renewcommand{\hat}{\widehat}
\renewcommand{\bar}{\overline}
\renewcommand{\epsilon}{\varepsilon}
\renewcommand{\leq}{\leqslant}
\renewcommand{\geq}{\geqslant}
\renewcommand{\phi}{\varphi}

\renewcommand{\Re}{\EuScript{R}e}

\newcommand{\lps}{\langle \! \langle}
\newcommand{\rps}{\rangle \! \rangle}

\begin{document}
\title[An optimal inequality]{An optimal inequality between scalar curvature and spectrum of
the Laplacian}
\author[H. Davaux]{H\'el\`ene Davaux} 
\date{\today}
\address{D\'epartement de Math\'ematique cc51, Universit\'e Montpellier II,  Laboratoire GTA, CNRS UMR 5030, F-34095 Montpellier, France}
\email{davaux@math.univ-montp2.fr}
\thanks{ }
\keywords{Dirac operator, scalar curvature, Kato inequality, 
index Theorem, spectral flow,
 eta invariant, Von Neuman algebra, spectrum of the Laplacian.
}
\subjclass{[2000] Primary:  58J50 ; Secondary: 35P15, 46L10, 58G11}
\begin{abstract}
For a Riemannian closed spin manifold and
under some topological assumption (non-zero $\widehat{A}$-genus or 
enlargeability in the sense of Gromov-Lawson), we give an optimal 
upper bound for  the infimum
of the scalar curvature in terms of the first eigenvalue of the Laplacian.
The main difficulty lies in the study of the odd-dimensional case.
On the other hand, we study the equality case for the closed spin Riemannian 
manifolds with non-zero $\widehat{A}$-genus.
This work improves an inequality which was first proved by
K. Ono in $1988$.
\end{abstract}

%
\maketitle
\setcounter{tocdepth}{1}
%
%
\section{Introduction}

Our result is in the framework of spin geometry, in which Dirac
operators and index Theorems play an important role (see \cite{LM89}).

Recall that, in the sixties, by applying the index Theorem
to the Dirac operator, A. Lichnerowicz \cite{Lic63} 
found a striking topological
obstruction to the existence of Riemannian metrics
of positive scalar curvature on a closed spin manifold.
Namely, if the $\hat A$-genus is not zero, such metrics do not exist.

In spite of its beauty, this result was somewhat frustrating:
it works in dimension $4k$ only, and does not give
anything for such  manifolds as tori, whose
$\hat A$-genus vanishes. The second break-through,
was made by M. Gromov and H.B. Lawson in the late seventies (see \cite{GL80}).
Their idea was to write down a Bochner-Weitzenb\"ock formula,
not for the classical Dirac operator
as Lichnerowicz did, but
for some <<twisted>> operator, defined
on the tensor product of the spin bundle by a suitable
vector bundle. Then, when there exist non homologically
trivial such bundles with arbitrarily small curvature,
they obtained an obstruction to the existence
of positive scalar curvature metrics.

More specifically, their obstruction works for
\emph{enlargeable manifolds}:

\begin{Edefinition} 
\label{Enlargeable}
Let $(V,g)$ be an $n$-dimensional closed oriented Riemannian manifold.
The manifold $V$ is said to be \emph{enlargeable} if for each constant 
$\epsilon >0$, there
exists a finite Riemannian spin covering $\tV_\epsilon$ of $V$ such that
$\tV_\epsilon$ admits an $\epsilon$-contracting map $f_\epsilon : 
\tV_\epsilon \rightarrow S^n$ (that is $\| {f_\epsilon}_* v \| \leq 
\epsilon \|v\|$ for all tangent vectors $v$ on $\tV_\epsilon$)
of non-zero degree.
\end{Edefinition}

Then for a closed spin manifold with non-zero $\widehat{A}$-genus
(A. Lichnerowicz) or for a closed  enlargeable  manifold
(Gromov-Lawson), we have the  inequality
\begin{equation*}
\label{Inequality}
\inf \scal(V,g) \leq 0 \textrm{ for any metric}
\end{equation*}
where $\scal$ denotes the scalar curvature.

In \cite{Ono88}, K. Ono improved this inequality.
Let us first introduce   
the total $\Ah$-class, $\Ah(V)$, of a $n$-dimensional 
manifold $V$ (\cite[pp. 231-233]{LM89}). 
Then the $\Ah$-genus, denoted by $\Ah_V[V]$, is just the evaluation 
of the class 
$\Ah(V) \in \H^{4*}(V, \QQ)$ on the 
fundamental homology class $[V] \in \H_n(V)$ of the manifold. 
The greatest lower bound of  the spectrum 
of the positive  Laplacian will be denoted by $\lambda_0$.

\begin{theorem}[K. Ono] 
\label{Ono}
Let $(V,g)$ be a closed Riemannian  manifold and 
$(\widetilde{V}, \widetilde{g})$ be its universal covering. 

\noindent
\textbf{(a)} 
If $V$ is spin with non-zero $\widehat{A}$-genus  
then
\begin{equation*}
\label{OnoInequality}
\inf \scal (V,g) \leq -4  \lambda_0(\widetilde{V}, 
\widetilde{g}).
\end{equation*}

\noindent
\textbf{(b)}
If $V$ is enlargeable, then the same 
inequality holds.
\end{theorem}

K. Ono relates his result with the following result 
of R. Brooks \cite{Bro81}
which says that for a closed Riemannian manifold $V$, 
$\lambda_0(\tV,\widetilde{g})=0$ if and only if the fundamental 
group of $V$ is amenable.
Here we call a group $G$ amenable if it has a left invariant mean.
Since
$\lambda_0(\tV, \tilde{g})$ is always non-negative, the result is 
particularly interesting
for a manifold of  non amenable fundamental group. 
Indeed a closed spin manifold with non-amenable fundamental group which
admits a metric of zero scalar curvature will have necessary zero 
$\widehat{A}$-genus.

In this paper we improve both statements {\it \textbf{(a)}} and 
{\it \textbf{(b)}} in  
 Theorem \ref{Ono} to optimal inequalities.

\begin{theo}
\label{AgenusCase} 
Let $(V,g)$ be a closed Riemannian spin manifold of dimension $n=4k$ and  
$(\widetilde{V}, \widetilde{g})$ be its universal covering. 
If $V$ satisfies $\Ah_V[V] \neq 0$, then 
\begin{equation}
\label{TheInequality}
 \inf \scal (V,g) \leq -4 \frac{n}{n-1} \lambda_0(\widetilde{V}, 
\widetilde{g}).
\end{equation}
\end{theo}

\begin{theo}
\label{EnlargeableCase} 
Let $(V,g)$ be a closed Riemannian manifold of dimension $n$  and  
$(\widetilde{V}, \widetilde{g})$ be its universal covering.
If $V$ is enlargeable, then the same 
inequality holds.
\end{theo}

While K. Ono's proof is easily adapted in the even dimensional
case (we mainly use the refined  Kato inequality instead of the 
classical one), the odd dimensional case requires more work, to which 
is devoted
the main part of this paper. 
We point out that, in the odd dimensional case,  K. Ono considered 
$V \times S^1$ to bring the problem back to even dimension. 
But in our case, this method gives rise to  the 
constant $(n+1)/n$ instead of $n/(n-1)$.

Now, let us make some remarks on the equality case.
Since the bottom of the spectrum of the hyperbolic space equipped 
with the canonical metric is just 
$(n-1)^2/4$, the  inequality
(\ref{TheInequality})  is optimal in the case of closed hyperbolic
(enlargeable) manifolds. Furthermore, under the assumptions of 
Theorem \ref{AgenusCase}, we can make  a careful study of 
the \emph{equality case}:

\begin{theo}
\label{EqualityCase}
Let $(V,g)$ be a closed Riemannian spin manifold of dimension $n=4k$, with 
non-zero $\widehat{A}$-genus. We have
$$\inf \scal(V,g) = -4 \frac{n}{n-1} \lambda_0(\tV, \tilde{g}) 
\Longleftrightarrow
\scal (V,g) = 0.$$
\end{theo}

Recall that A.Futaki showed in \cite{Fut93} that for
a manifold $(V,g)$ which satisfies the assumptions of Theorem \ref{AgenusCase},
if $\# \pi_1(V) |\Ah_V[V]| > 2^{n/4}$ then $V$ does not admit any metric with
$\scal(V,g) \geq 0$.
Therefore in the equality case of  Theorem \ref{AgenusCase} we have necessarily
  $\# \pi_1(V) |\Ah_V[V]| \leq 2^{n/4}$.

When $\pi_1(V)$ is finite 
(and therefore when $\lambda_0(\tV, \tilde{g})$ is zero) 
we can use $\lambda_1(\tV, \tilde{g})$ instead of $\lambda_0(\tV, \tilde{g})$.
This improvement is due to V. Mathai in \cite[Theo. 2.5]{Mat92} who was
inspired from Meyer's Lemma revisited by P. B\'erard in 
\cite[Appendix III]{Ber88}.

\section{Outline of the proofs}
\label{Outline}

Theorem \ref{Ono} is originally due to K. Ono \cite{Ono88};
however we prefer a shorter proof given by V. Mathai in
\cite[Theo. 2.5]{Mat92} and it is his proof which we shall update.

First we introduce some notation: $V$ will always be a spin closed Riemannian 
$n$-dimensional manifold and $(\tV, \tilde{g})$ a  Riemannian covering
(the universal Riemaniann covering unless otherwise stated). In the following
the tilde  is always used to denote the lift of the object
under consideration to the covering $\tV$. Since $V$ is spin we can consider
the spinor bundle $\Spin$ over $V$ and the associated Dirac operator $\Dirac$
acting on spinors over $V$. If  $Y$ is a Riemannian bundle over $V$
with a Riemannian connection,
we can construct the Clifford bundle $\Spin \otimes Y$ endowed with
the associated twisted Dirac operator $\Dirac^Y$ acting on sections
of the bundle $\Spin \otimes Y$.

In the present article we will prove the inequality
\begin{equation*}
\inf \scal (V,g) \leq -4 \frac{n}{n-1} \lambda_0(\tV, \tilde{g})
\end{equation*}
with various assumptions on $V$. But the proofs are essentially similar and
can be
decomposed into three  key points:
\begin{enumerate}
\item \textbf{We find a condition} under which there exists
a Dirac operator $\tD$ on $\tV$ 
such that $0$ belongs to the spectrum of this operator. 
In order to use the action of $\pi_1(V)$,
we look for such an operator as a lifted operator of 
a Dirac operator $\D$ on $V$.

In the even
dimensional case
we shall use  
Atiyah's $\G$-index Theorem  (see \cite{Ati76} and \cite[chap. 15]{Roe98}). 
For the odd dimensional case we need the spectral flow
introduced in \cite[$\S 7$]{APS76} 
(see also \cite[$\S 6 \frac{11}{12}$]{Gro96} and 
\cite[$\S 1$]{Nic99}) and the $\G$-index theorem developed by M. Ramachandran
in \cite{Ram93} for the non local Atiyah-Patodi-Singer boundary problem on
 non compact manifold with boundary.

\item  \textbf{From this, we prove an approximate inequality}. 
The Dirac operators under 
consideration are Dirac operators  twisted by a bundle denoted $Y$. 
Such operators give rise to a Bochner-Weitzenb\"ock formula 
(see \cite[$\S2$]{GL83}). Then, 
the refined Kato inequality \cite{CGH00} combined with a standard  
Rayleigh quotient argument gives the following inequality, 
which is nearly what we want: 
$$\inf \scal (V,g) \leq -4 \frac{n}{n-1} \lambda_0(\tV, \tilde{g}) + \alpha_n
\|R^Y\|$$
where $\alpha_n$ is some universal constant and $R^Y$ 
is the curvature of the bundle $Y$.

\item \textbf{We construct},
for any $\epsilon >0$, 
a sequence of bundles $Y_\epsilon$ on $V$
  such that
$0$ belongs to the spectrum of $\tilde{\Dirac^{Y_\epsilon}}$ 
(i.e. satisfies the
condition  found previously in the first step (1))
and $\|R^{Y_\epsilon}\| \leq \epsilon$. 
Letting $\epsilon \rightarrow 0$, we obtain the desired inequality 
(\ref{TheInequality}).
Here the enlargeability of $V$ is crucial.
\end{enumerate}

\vspace{0.3cm}
The paper is organized as follows.
In Section \ref{GIndex} we recall 
the $\G$-index theory and achieve the first step of the proof, especially
in the odd dimensional case.
Section \ref{RKato} is devoted to the second
point (which does not depend on the dimension).
In Section \ref{Proof1} we prove Theorem \ref{AgenusCase} 
and its 
equality case, Theorem \ref{EqualityCase}. 
In Section \ref{Proof2} we give the proof of Theorem 
\ref{EnlargeableCase} in the
even and odd dimensional cases.
In Section \ref{Generalization} we make some remarks and  generalizations 
of our results.

\section{Use of $\G$-index}
\label{GIndex}

\subsection{Atiyah's $\G$-index Theorem}

For this section we refer to the original
article of M. Atiyah \cite{Ati76} (see also
 \cite[chap. $15$]{Roe98}). 
The general context is the following:
$\tV$ denotes a Galois covering of $V$ with Galois group $\G$.
 The group $\G$ acts discontinuously on $\tV$
by deck transformations and $\tV / \G=V$. If $S$ is a Clifford bundle over
$V$ with Dirac operator $\D$, let $\tilde{S}$ and $\tD$ denote their 
natural lifts on $\tV$.
In the following $\G$ is almost always the fundamental
group of $V$.

On the \emph{a priori} non-compact manifold $\tV$, 
a central idea is to introduce 
an algebra of smoothing operators on $\tV$ that reflects the extra 
structure of the $\G$-action.

\begin{Edefinition}
\label{DefiHS}
With notation as above, we define the set of 
\emph{$\G$-Hilbert-Schmidt operators} $A$ on $\L^2(\widetilde{S})$
by the following conditions:
\begin{enumerate}
\item[($i$)] $A$ is bounded on $\L^2(\widetilde{S})$.
\item[($ii$)] $A$ is $\G$-invariant; that is, for all $s \in 
\L^2(\widetilde{S}), \; A(\gamma s)=\gamma(As)$ where by definition
$\gamma s(\tilde{x})=s(\gamma^{-1}\tilde{x})$.
\item[($iii$)] $A$ is represented by a smooth kernel $k(\tilde{x},\tilde{y})$
so that $As(\tilde{x})= \int_{\tV} k(\tilde{x},\tilde{y}) s(\tilde{y}) 
\d \tilde{y}$. Moreover $\tilde{x} \mapsto k(\tilde{x}, \cdot)$ and
$\tilde{y} \mapsto k(\cdot,\tilde{y})$ are smooth maps of
$\tV$ into the Hilbert space $\L^2(\widetilde{S})$.
\item[($iv$)] There is a constant $C$ such that for all $\tilde x \in \tV$, 
we have
$\int_{\tV} |k(\tilde{x},\tilde{y})|^2\d \tilde{y} <C$ and  
for all $\tilde y \in \tV$, 
we have $\int_{\tV} 
|k(\tilde{x},\tilde{y})|^2 \d \tilde{x} < C$.
\end{enumerate}
\end{Edefinition}

By the Riesz representation Theorem for functionals
on Hilbert space, the quantity $\int_{\tV} |k(\tilde{x},\tilde{y})|^2
\d \tilde{y}$ 
is the square
of the norm of the linear functional $s \mapsto As(\tilde{x})$ 
on $\L^2(\widetilde{S})$. Therefore, there is a constant $C$ such that
$\int_{\tV} |k(\tilde{x},\tilde{y})|^2\d \tilde{y} <C$ if and only if $A$ maps
 $\L^2(\widetilde{S})$
continuously to the space $\CB(\widetilde{S})$ of bounded continuous
sections of $\widetilde{S}$. Similarly there is a constant $C$ such
that $\int_{\tV} |k(\tilde{x},\tilde{y})|^2\d \tilde{x}<C$ if and
only if $A^*$ maps  $\L^2(\widetilde{S})$ continuously to $\CB(\widetilde{S})$.

The set of $\G$-Hilbert-Schmidt operators is an algebra. 

An operator $A$ is said to be \emph{of $\G$-trace} if there exist
 two $\G$-Hilbert-Schmidt operators $B_1$ and $B_2$ such that
$A=B_1B_2$.
If $A$ is an operator of $\G$-trace, let $k$ be its kernel and choose
any fundamental domain $F$ for the $\G$-action on $\tV$; then
$$\Tr_\G(A) := \int_F \Tr k(\tilde{x},\tilde{x}) \d \vol (\tilde{x})<\infty.$$

Now we give an important property of the $\G$-trace 
(nearly the same as for the standard trace of matrix): 
for  $A$ and $B$ two $\G$-Hilbert-Schmidt operators or $A$ an operator of 
$\G$-trace and $B$ a bounded operator, the operators $AB$ and 
$BA$ are of $\G$-trace 
and
$$\Tr_\G (AB)= \Tr_\G (BA).$$

We also need to  recall the \emph{$\G$-dimension} which 
is motivated by the fact that the
trace of a projection operator is the dimension of its image. 
If $H$ is a subspace of $\L^2(\tilde{S})$ with the property that the 
orthogonal projection operator $P$ from $\L^2(\tilde{S})$ onto
$H$ is $\G$-Hilbert-Schmidt, then we define
$$\dim_\G(H) := \Tr_\G(P).$$

Finally we restate the \emph{Atiyah's $\G$-index Theorem} \cite{Ati76}:
suppose $\D$ (and so $\tD$) is a $\ZZ_2$-graded Dirac operator. 
The orthogonal projection $P$ onto the kernel of $\tD$ is
$\G$-Hilbert-Schmidt and
$$\ind_\G(\tD^+) := \dim_\G(\ker \tD^+) - \dim_\G(\ker \tD^-)=\ind(\D^+).$$

We recall that a $\ZZ_2$-graded operator $\D$ from $\L^2(S)$ to $\L^2(S)$ is
an operator which decomposes itself as 
$$\left( \begin{array}{cc} 0 &  \D^-\\ \D^+ &0\end{array}\right)$$
on the eigenspaces $S^+$ and $S^-$ of an involution $\epsilon$ acting on $S$.
For example, the standard Dirac operator on a spin even-dimensional
manifold $V^{2m}$ is $\ZZ_2$-graded. 
Here $\epsilon$ is the Clifford multiplication by the volume
form $\omega = i^{m} e_1 \cdot e_2 \cdots e_{2m}$ (where
the $(e_i)_{i=1,\cdots,2m}$ form an orthonormal basis of $T^*V$) acting on the 
spinor bundle.

\subsection{First step in the even dimensional case}

In the even dimensional case we can now give a condition which ensures
the existence of a Dirac operator on $\tV$ such that $0$ belongs to
its spectrum.

\begin{proposition}
\label{EvenCondition}
Let $\D$ be a $\ZZ_2$-graded Dirac operator on $V$ and $\tD$ be its lift.
If $\ind \D^+ \neq 0$, then $0$ belongs to the point spectrum of $\tD$. 
\end{proposition}
 
\begin{proof}
It is just an application of Atiyah's $\G$-index Theorem.
Indeed $\ind_{\G} \tD^+ \neq 0$ implies that $\dim_\G \ker(\tD) \neq 0$.
Using the definition of $\G$-dimension, this ensures that
there exists a non-zero spinor in this kernel.  
\end{proof}

\begin{remark}
\label{DiracCondition}
For the Dirac operator $\Dirac$ on spinors the condition is 
$\widehat{A}_V[V] \neq 0$, which is exactly the assumption of 
Theorem \ref{AgenusCase}. 
For the twisted Dirac operator $\Dirac^Y$, the condition
is $\{\widehat{A}(V) \ch(Y)\}[V] \neq 0$ (see \cite[Theo. $13.10$]{LM89}). 
\end{remark}

\subsection{First step in the odd dimensional case}
\label{Index}

On an odd dimensional manifold, the ordinary index of every 
elliptic operator  is always zero (\cite[Theo. $13.12$]{LM89}). 
In order to overcome this difficulty, 
we consider the spectral flow of a family of Dirac operators instead of 
the index Theorem of only one Dirac operator.  
The use of the spectral flow is inspired by the proof of Vafa-Witten's
Theorem on
eigenvalues of the Dirac operator, given in 
\cite{Ati85}.

\subsubsection{Condition on the spectral flow}

We start with a smooth family 
$\{\D_u\}_{u\in [0,1]}$ of Dirac operators on $V$ and lift it to $\tV$
 to obtain the family $\{\tD_u\}_{u\in [0,1]}$ of Dirac operators
on $\tV$. We search for a condition which ensures that there exists
some $u_0 \in [0,1]$ such that $0$ belongs to the spectrum of $\tD_{u_0}$. For
this we will use the adequate index Theorems on $V \times [0,1]$ and
$\tV \times [0,1]$. For $V \times [0,1]$,  we need the index Theorem of
Atiyah-Patodi-Singer (see \cite{APS76}) for a compact manifold with 
boundary and its associated $\eta$-invariant. We will  also recall what
is the spectral flow of a family of Dirac operators and
in particular, the link between
the $\eta$-invariant and the spectral flow.
For $\tV \times [0,1]$,
 we will use the $\G$-index Theorem developed by M. Ramachandran
in \cite{Ram93} for a non-compact manifold with boundary where the 
$\eta_{\G}$-invariant plays an important role. We will then relate
the derivatives of $\eta$ and $\eta_{\G}$ 
(\cite[$\S 1.7 \& 1.10$]{Gil84}). 
In the following we will denote the spectrum of an operator $\D$ by
$\spec(\D)$.

Our main task in this section is to show:

\begin{theorem}
\label{OddCondition}
If the spectral flow  of the smooth family 
$\{\D_u\}_{u\in [0,1]}$ of Dirac operators 
on $V$ is non-zero, then there exists $u_0 \in [0,1]$
such that $0$ belongs to the spectrum of $\tD_{u_0}$.
\end{theorem}

\begin{proof}
We borrow  here the presentation of the Atiyah-Patodi-Singer index Theorem 
and of the $\eta$-invariant from
L.I. Nicolaescu \cite[$\S 1$]{Nic99}.
On $V \times [0,1]$, we consider the operator $\mathcal{D} = \pdd{u} - \D_u$. 
The index Theorem with the spectral boundary condition of 
Atiyah-Patodi-Singer (\cite{APS76}) 
for the manifold $V \times [0,u]$
with boundary $V \times \{0\} \cup V \times \{u\}$, where 
$V \times \{0\}$ correspond
to $V$ negatively oriented and $V \times \{u\}$ correspond
to $V$ positively oriented (the orientation of $\partial(V \times [0,1])$ is 
given by the outer normal), 
 gives
\begin{eqnarray*}
\ind(\mathcal{D}, V\times [0,u]) = A(u) -\frac{1}{2} \left(
-\eta(0) + h(0) + \eta(u) + h(u) \right) 
\end{eqnarray*}
where $\eta(u)$ is the $\eta$-invariant of the manifold $(V, \D_u)$, $h(u)
= \dim \ker \D_u$ and $A(u)$ is the integral over $M \times [0,u]$
of the index density determined by the operator $\mathcal D$ which is 
a completely local object.

We introduce the $\xi$-invariant (or reduced $\eta$-invariant) by $\xi(u)=
(h(u)+\eta(u))/2$. Hence 
\begin{eqnarray}
\label{index1}
\ind(\mathcal{D}, V\times [0,u]) = A(u) - \left(
\xi(u) - \xi(0) \right) - h(0).
\end{eqnarray}

Recall what the \emph{spectral flow} is. 
The discontinuities of $\xi(u)$ are due to jumps in $h(u)$. We describe how
the jump in $h(u)$ affects $\xi(u)$ in a simple, yet generic situation. We
assume $\D_u$ is a regular family i.e.
\begin{itemize}
\item the resonance set $\mathcal{Z} = \{ u \in [0,1]; \; h(u) \neq 0 \}$ is
finite,
\item For every $u_0 \in \mathcal{Z}$, there exists $\epsilon >0$, an open 
neighbourhood $\mathcal{N}$ of $u_0$ in $[0,1]$ and smooth maps
$\lambda_k: \mathcal{N} \rightarrow ]-\epsilon, \epsilon[, \; k=1,\dots,h(u)$
such that for all $u \in \mathcal{N}$ the family $\{\lambda_k(u)\}_k$
describes \emph{all} the eigenvalues of $\D_u$ in $]-\epsilon,\epsilon[$
(including multiplicities) and, moreover, $\dd{u} \lambda_k(u_0) \neq 0$
for all $k=1,\dots,h(u_0)$.
\end{itemize}
Now for each $u \in \mathcal{Z}$, set $\sigma_\pm = \# \{k | \pm 
\dd{u} \lambda_k(u) >0 \}$ and
$$\Delta_u \sigma = \left\{ \begin{array}{l} 
-\sigma_-(0) \textrm{ if } u=0\\
\sigma_+(u)-\sigma_-(u) \textrm{ if } u\in ]0,1[\\
\sigma_+(1) \textrm{ if } u=1\\
\end{array} \right.$$
If $\Delta_u \xi : = \lim_{\epsilon \rightarrow 0^+}\left(\xi(u+\epsilon)
-\xi(u-\epsilon)\right)$, we see that $\Delta_u\xi =0$ if
$u \notin \mathcal{Z}$ while for $u \in \mathcal{Z}$ we have
$\Delta_u \xi = \Delta_u \sigma$. Finally, define the spectral flow
of the family $\D_u$ by 
\begin{eqnarray}
\label{sflowEq}
\sflow(1)=\sflow(\D_u, \, u \in [0,1]) = \sum_{u\in [0,1]} \Delta_u \sigma
=\sum_{u\in [0,1]} \Delta_u \xi.
\end{eqnarray}
In fact the spectral flow counts how many eigenvalues of $\D_u$ cross 
zero as $u$ moves from $0$ to $1$.

We rewrite the equality (\ref{index1}) as follows, 
\begin{eqnarray*}
\xi(u) - \xi(0) =A(u) -
\ind(\mathcal{D}, V\times [0,u])  -h(0).
\end{eqnarray*}
The term $A(u)$ depends smoothly on $u$. The term $\left(- 
\ind(\mathcal{D}, V\times [0,u])-h(0)\right)$ is  integer-valued so
it cannot be smooth, unless it is constant. If $\bar{\xi}(u) :=
\xi(u) \mod \ZZ$ (seen as a map from $[0,1]$ to  $[0,1]$ identified with 
$S^1$) 
then $u \mapsto
\bar{\xi}(u)$ is smooth. Therefore we deduct
\begin{equation}
\label{index2}
\sflow(u) 
\stackrel{(\ref{sflowEq})}{=}\sum_{v\in [0,u]} \Delta_v \xi = - 
\ind(\mathcal{D}, M \times [0,u])-h(u)  
\end{equation}
and
\begin{equation}
\label{index3}
\dd{u} \bar{\xi}(u) = \dd{u} A(u) \textrm{ and so } A(u) = \int_0^u 
\dd{v} \bar{\xi}(v) \d v. 
\end{equation}

Moreover we assume $\D_0 = \D_1$, hence
{\setlength\arraycolsep{2pt}
\begin{eqnarray*}
\ind(\mathcal{D}, V\times [0,1]) & \stackrel{(\ref{index1})}{=} & 
A(1) - \frac{1}{2}(h(0)+h(1)) \stackrel{(\textrm{hyp.})}{=} 
A(1)- h(0) 
\stackrel{(\ref{index3})}{=} \int_0^1 
\dd{u} \bar{\xi}(u) \d u 
- h(0)\\
 & \stackrel{(\ref{index2})}{=} &  - \sflow(1) - h(0). 
\end{eqnarray*}}
Thus we have proved the following Lemma:
\begin{lemma}
\label{sflow}
With the previous notations:
 $$A(1) = -\sflow(1) 
=\int_0^1 \dd{u} \bar{\xi}(u) \d u = \frac{1}{2} 
\int_0^1\dd u \bar{\eta}(u) \d u$$
where $\bar{\eta}(u)= 2 \bar{\xi}(u)$ (representant of $\eta(u)$ in $[0,2[$
-- seen as a map from $[0,1]$ to $S^1$).
\end{lemma}

\begin{remark}
\label{Convention}
The minus sign  in the equality $A(1) = -\sflow(1)$ is due to the fact that 
we use $\mathcal{D}= \pdd{u} - \D_u$ instead of 
the convention of Atiyah-Patodi-Singer \cite{APS76} which is to use
$\mathcal{D}= \pdd{u} + \D_u$. Moreover we do not use the same orientation 
on $V \times [0,1]$ as  in \cite{APS76}. This different convention
on orientation is also used by Bismut-Freed 
\cite[Theo. $2.11$, Rem. $6$]{BF86}.
\end{remark}

We now study
$\tilde{\mathcal{D}}= \pdd{u}- \tD_u$ on $\tV \times [0,1]$. 
The relevant index 
Theorem in this case gives  the formula \cite{Ram93}
$$\ind_\G(\tilde{\mathcal{D}}, \tV \times [0,u]) = A(u) -
\frac{1}{2} \left(-\eta_\G(0) + h_\G(0) + \eta_\G(u) + h_\G(u) \right)$$
where $\eta_\G(u)$ is the $\eta_\G$-invariant of the manifold $(\tV, \tD_u)$
(defined in the next Sub-section \ref{SectionEta}),
 $h_\G(u)= \dim_\G \ker \tD_u$ and $A(u)$ is the same as 
in the formula (\ref{index1}) (in particular it is an integral over $V 
\times [0,u]$ and not on $\tV \times [0,u]$).

By contradiction,
we assume that for all $u \in [0,1]$, 
$0 \notin \spec \D_u$
and thus $h_\G(u)=0$. 
In the following Sub-section \ref{SectionEta}, we shall show the
crucial lemma:

\begin{lemma}
\label{Eta}
Assume $0 \notin \spec(\tD_u)$ for all $u$. Then
the
$\eta_\G$-invariant is a smooth function of $u$ satisfying 
$$\dd{u} \eta_\G(u)= 
\dd{u} \bar{\eta}(u).$$
\end{lemma}

Hence, as $\D_0=\D_1$ (and so $\tD_0=\tD_1$), using Lemma \ref{sflow}
and Lemma \ref{Eta}, 
we have
$$0=\eta_\G(1)-\eta_\G(0)= \int_0^1 \!\! \dd{u} \eta_\G(u) \d u = 
\int_0^1 \!\! 
\dd{u} \bar{\eta}(u) \d u = - 2 \, \sflow(1).$$
Therefore we have found a contradiction and Theorem \ref{OddCondition} 
is proved.  
\end{proof} 

\begin{remark}
\label{Bismut-Freed}
We emphasize that the spectral flow of a family of Dirac operators
does not depend on the index Theorem of the compact  manifold with boundary
$V \times [0,1]$ but on the classical Atiyah-Singer index Theorem 
on the closed manifold $V \times S^1$ (see \cite[Theo. $2.11$]{BF86}). In fact,
$\sflow(\D_u; u\in [0,1], D_0=D_1) = -\ind(\mathcal{D})$ for any
operator on $V \times [0,1]$ such that $\mathcal{D}|_{V \times \{u\}}
=-\D_u$. This remark will be important when we generalize our result 
to the infinite $\Karea$ case in Section \ref{Generalization}.
\end{remark}

\subsubsection{Study of the $\eta_{(\G)}$-invariant}%
\label{SectionEta}

The proof of Lemma \ref{Eta} will occupy the whole of this section.
It is based on 
the same principle as  the proof of Atiyah's $\G$-index
Theorem given in \cite[Chap. $15$]{Roe98}. He used an asymptotic expansion
 of the kernel of $e^{-t\tD^2}$ to prove that $\ind_\G \tD^+$ is 
given by an integral on a fundamental domain
only involving the symbol of $\tD$ and therefore he concluded to the index
equality. Therefore, here
our main task is to show that, in our case,  
 $\dd{u} \eta_\G(u)$ is a local invariant which can be computed as 
an integral of a local form. 

\begin{lemma}
\label{EtaLocal}
When for every $u \in [0,1]$, $0 \notin \spec(\tD_u)$, 
$$\dd{u} \eta_\G(u)= \int_F \! \! a_n(x, \dd{u} \tD_u, \tD_u) \d x$$
where $F$ is a fundamental domain and $a_n$ is a local invariant in the jets
of the symbol of $(\dd{u} \tD_u,\tD_u)$. 
\end{lemma}

The result for $\bar{\eta}$ is already known 
\cite[Lem. $1.10.3$]{Gil84} and corresponds to the case $\G=\{id\}$.
As the symbol of the operators and their lifts
are the same and as we integrate over a fundamental domain, we conclude that 
Lemma \ref{EtaLocal} implies Lemma \ref{Eta}.

\begin{DefProp}[M. Ramachandran \cite{Ram93} $\S 3.1$]
\label{DefEta}
The $\eta_\G$-invariant is defined by
$$\eta_\G(u)= \frac{1}{\G(1/2)} \int_0^\infty \!\! t^{-1/2} \Tr_\G(\tD_u 
e^{-t \tD^2_u} )\d t.$$
That is, 
the following limits
$$\frac{1}{\G(1/2)} \lim_{T \rightarrow \infty}  \lim_{\epsilon \rightarrow 0}
 \int_\epsilon^T t^{-1/2} \Tr_\G(\tD_u 
e^{-t \tD^2_u} )\d t $$
are finite. Moreover the limits are uniform with respect to $u$.
\end{DefProp}

\begin{proof}
For $0$ the convergence is based  on a  Bismut-Freed estimation 
\cite[Theo. $2.4$]{BF86}:
denoting by  $K(t,x,y)$ the 
kernel of the integral operator $e^{-t \tD^2}$  then
$$ | \Tr_x \tD K(t,x,x) | < C t^{1/2}$$
where $C$ is a constant depending on the local geometry of $V, S$, the
dimension of $V$ and the rank of $S$. Hence the limit for 
$\epsilon \rightarrow 0$ 
is well defined. Moreover,
as $C$ does not depend on $u$, the convergence 
is uniform with respect to $u$. 

For $+\infty$ we use the same standard spectral measure argument 
as M. Ramachandran 
in \cite[Theo. $3.1.1$]{Ram93}. The convergence 
is again uniform with 
respect to $u$.

As a conclusion,  the $\eta_\G$-invariant is well defined and continuous.
\end{proof}

Now we are ready to prove Lemma \ref{EtaLocal}.

\begin{proof} 
The proof is decomposed in three steps from ($a$) to ($c$).
\emph{A priori}, we do not know that $\dd{u} \eta_\G(u)$ exists. Its existence
will be a consequence of the proof.

($a$) We begin with establishing the formula:
\begin{equation}
\label{Formula}
\dd{u} \Tr_\G(\tD_u e^{-t \tD^2_u}) = (1+2 t \dd{t}) \Tr_\G (\dd{u} 
(\tD_u) e^{-t \tD^2_u}).
\end{equation}

The proof of this formula uses especially
the article \cite[\S $4$]{CG85}.

Before starting, we make one remark. Subsequently, all inversions
of limits and integrals, integrals and integrals, derivations and integrals, 
derivations and limits are permitted by the fact that we work on compacts
($u \in [0,1]$, the $\G$-trace is an integral on a compact fundamental
domain $F$,\dots) and the fact that the functions under consideration
 are smooth (in all variables).

Now we can start by writing Duhamel's principle
\begin{eqnarray}
\label{1}
\lefteqn{ \tD_0 e^{-(t-\epsilon) \tD_u^2}e^{-\epsilon \tD_0^2} -\tD_0 e^{-\epsilon \tD_u^2}e^{-(t-\epsilon) \tD_0^2}} \nonumber \\
& & {}  =  -\int_\epsilon^{t-\epsilon} \!\!
\dd{s} \left[ \tD_0 e^{-(t-s) \tD_u^2}e^{-s \tD_0^2} \right] \d s  \nonumber \\
& & {} = - \int_\epsilon^{t-\epsilon} \!\! 
\left[ \tD_0 \tD_u^2 e^{-(t-s)\tD_u^2}
e^{-s\tD^2_0} -   \tD_0 e^{-(t-s)\tD_u^2}
e^{-s\tD^2_0} \tD^2_0 \right] \d s. 
\end{eqnarray}
If we take the $\G$-trace, the second term of the last line can be rewritten
\begin{eqnarray}
\label{2}
\lefteqn{ \int_\epsilon^{t-\epsilon} \!\! \Tr_\G \left( \tD_0 e^{-(t-s)\tD_u^2}
e^{-s\tD^2_0} \tD^2_0 \right) \d s } \nonumber\\
& & {}  = \int_\epsilon^{t-\epsilon}\!\!\Tr_\G \left( (\tD_0 e^{-(t-s)\tD_u^2}
e^{-\frac{s}{2}\tD^2_0})(e^{-\frac{s}{2}\tD^2_0} \tD^2_0) \right) \d s
\nonumber\\
& & {}  = \int_\epsilon^{t-\epsilon} \!\!\left((e^{-\frac{s}{2}\tD^2_0} \tD^2_0)
(\tD_0 e^{-(t-s)\tD_u^2}e^{-\frac{s}{2}\tD^2_0})  \right) \d s\nonumber\\
& & {}  = \int_\epsilon^{t-\epsilon} \!\! \Tr_\G \left( \tD_0^3 e^{-(t-s)\tD_u^2}
e^{-s \tD_0^2} \right) \d s.
\end{eqnarray}
In the previous computation we have used the commutation 
property of the $\G$-trace but
we need to verify the assumptions,  that  is the two operators
inverted $(\tD_0 e^{-(t-s)\tD_u^2}
e^{-\frac{s}{2}\tD^2_0})$ and $(e^{-\frac{s}{2}\tD^2_0} \tD^2_0)$
 are $\G$-Hilbert-Schmidt. It is the aim of the following sub-lemma:

\begin{sublemma}
\label{bounded}
Let $S$ be a Clifford bundle over a compact Riemannian manifold $V$ 
of dimension $n$. Let $D :
C^\infty(S) \rightarrow C^\infty(S)$ be a  Dirac operator
and  $Q : C^\infty(S) \rightarrow C^\infty(S)$ be
an auxiliary partial differential operator of order $a \geq 0$.
We lift all these objects to a $\G$-covering $\tV$ ($\tV$ can be $V$ itself).

Then the operator
$\tilde{Q} e^{-t\tilde{D}^2}$ is a well defined (i.e. 
$\Im e^{-t\tilde{D}^2} \subset \dom \tilde{Q}$), infinitely
smoothing operator with $\G$-invariant kernel 
$k(t,\tilde{x},\tilde{y})=\tilde{Q}_x K(t,\tilde{x},\tilde{y})$ 
(where $K(t,\tilde{x},\tilde{y})$ is the kernel of $e^{-t\tilde{D}^2}$
and the operator $\tilde{Q}_x$
acts on  $K(t,\tilde{x},\tilde{y})$, 
considered as a function of $\tilde{x}$).  
Moreover  $\tilde{Q} e^{-t\tD^2}$ is a 
bounded operator from $L^2(\tilde{S})$ to $L^2(\tilde{S})$
and even a $\G$-Hilbert-Schmidt operator.
\end{sublemma}

\begin{proof}
We prove successively all the points of the sub-lemma. 

First we have 
$\Im e^{-t\tilde{D}^2}\subset\dom\tilde{D}^k$. Indeed
the functional calculus ensures that the operator $f_k(\tD)=
\tD^ke^{-t\tD^2}$ for all $k$ where $f_k(x)=x^ke^{-tx^2}$ are well defined 
because the functions $f_k$ are rapidly 
decreasing. 
On the other hand the non-compact Sobolev embedding Theorem 
and elliptic estimates
(see \cite[Prop. $15.4$]{Roe98} for a proof) assert
that for any $p >\frac{n}{2}$ and $r \geq 0$, there is a constant $c$ such
that 
$$\|s\|_{CB^r} \leq c (\|s\|+ \|\tD s\| +\cdots +\|\tD^{p+r}s\|)$$  
for all $s$ in the domain of $\tD^{p+r}$ where 
$\mathrm{CB}^r(\tilde{S})$ is the space of
sections of $\tilde{S}$ which are $r$ times continuously
differentiable with all derivatives bounded, $\|s\|_{CB^r}$ is
the sum of the supremum norms of all the derivatives of $s$ and
$\|.\|$ is the $\L^2$-norm.
This proves that
$\dom \tD^{p+r} \subset \mathrm{CB}^r(\tilde{S})$.
Thus we have proved the first point: $\tilde{Q} e^{-t \tilde{D}^2}$ is
well defined. 

Moreover,  for any rapidly 
decreasing function $f$, the functional calculus
(see \cite[Prop. $9.20$]{Roe98} for a proof)
estimates the norm of the operator $f(\tD)$ from $\L^2(\tilde{S})$ to
$\L^2(\tilde S)$ by
$ \displaystyle \|f(\tD)\| \leq \sup |f(\lambda)| \textrm{ for } \lambda \textrm{ in the 
spectrum of } \tD.$ 
Combining
this with  the non-compact Sobolev embedding Theorem gives that
 $e^{-t\tD^2}$ maps
$L^2(\tilde{S})$ continuously in $\mathrm{CB}^r(\tilde{S})$ for all $r$.
In addition,
$\tilde{Q}$ maps trivially $\mathrm{CB}^r(\tilde{S})$ continuously
in $\mathrm{CB}^{r-a}(\tilde{S})$. Thus for all $k$, $\tilde{Q} e^{-t\tD^2}$
maps $L^2(\tilde{S})$ continuously in $\mathrm{CB}^k(\tilde{S})$:
$\tilde{Q} e^{-t\tD^2}$ is smoothing. 

An easy computation gives the announced 
kernel and the fact that it is $\G$-invariant. 

It remains to show that $\tilde{Q} e^{-t\tD^2}: L^2(\tilde{S}) \rightarrow
L^2(\tilde{S})$ is bounded. For this we use the closed graph Theorem i.e. 
consider a sequence $(u_n, \tilde{Q} e^{-t\tD^2}u_n)$ in the graph
which converges to $(u,v)$ in $L^2$-norm, we must show that
$\tilde{Q} e^{-t\tD^2}u=v$. But we know that the sequence
$\tilde{Q} e^{-t\tD^2}u_n$ converges to $\tilde{Q} e^{-t\tD^2}u$ in 
$\mathrm{CB}^r$-norm and therefore converges pointwise. By assumptions
the sequence $\tilde{Q} e^{-t\tD^2}u_n$ converges to $v$ in $L^2$-norm 
and so there exists a (dominated) subsequence that converges almost everywhere.
By uniqueness of the  limit 
$\tilde{Q} e^{-t\tD^2}u=v$.

Reading again Definition \ref{DefiHS} of a $\G$-Hilbert-Schmidt operator 
and the paragraph that follows it, 
we conclude  
that $\tilde{Q} e^{-t\tD^2}$ is a $\G$-Hilbert-Schmidt operator.
\end{proof}

\noindent
By definition of $\nabla_u$, all the space $CB^r(\tilde{S}, 
\tilde{\nabla}_u)$ are equal because the norms are equivalent. Therefore,
 $\tilde Q=\tD_0 e^{-(t-s)\tD_u^2}$ is well defined and we can apply  
Sub-lemma \ref{bounded} to $\tilde Q e^{-\frac{s}{2}\tD^2_0} =
(\tD_0 e^{-(t-s)\tD_u^2})
e^{-\frac{s}{2}\tD^2_0}$ as well as, directly to $e^{-\frac{s}{2}\tD^2_0}
\tD^2_0 $
which is equal to 
 $\tilde Q e^{-\frac{s}{2}\tD^2_0} =\tD^2_0 e^{-\frac{s}{2}\tD^2_0}$.

Let return to the computation.
The equality $(\ref{1})$ yields
\begin{eqnarray}
\label{3}
\lefteqn{ \Tr_\G \left( \tD_0 e^{-(t-\epsilon) \tD_u^2} e^{-\epsilon \tD_0^2}
\right) - \Tr_\G \left(  \tD_0 e^{-\epsilon \tD_u^2} e^{-(t-\epsilon) \tD_0^2}
\right) } \nonumber\\
& & {}  =  - \int_\epsilon^{t-\epsilon} \!\! 
\left[ \Tr_\G \left( \tD_0 \tD_u^2 
e^{-(t-s)\tD_u^2}e^{-s\tD^2_0} \right) - \Tr_\G \left(  \tD_0 e^{-(t-s)\tD_u^2}
e^{-s\tD^2_0} \tD^2_0 \right) \right] \d s  
\nonumber\\
& & {}  = \int_\epsilon^{t-\epsilon} \!\! 
\Tr_\G \left( \tD_0 (\tD_0^2 -\tD_u^2) 
e^{-(t-s)\tD_u^2}e^{-s\tD^2_0} \right) \d s.
\end{eqnarray}
If we differentiate (\ref{3}) with respect to $u$ and set $u=0$, the 
right-hand side becomes
\begin{eqnarray*}
\label{4}
 -  \int_\epsilon^{t-\epsilon}\!\!\Tr_\G \left( \tD_0 \dd{u}(\tD_u^2) 
e^{-(t-s)\tD_u^2}e^{-s\tD^2_0} \right) \d s &=& -  \int_\epsilon^{t-\epsilon}\!\!\Tr_\G \left( \tD_0 \dd{u}(\tD_u^2) 
e^{-t\tD_u^2} \right) \d s 
\end{eqnarray*}
\begin{eqnarray*}
& & {}  = -(t-2\epsilon)  \Tr_\G \left(  \tD_0 \dd{u}(\tD_u^2)_{|u=0} 
e^{-t\tD_0^2} \right)\nonumber\\
& & {}  = -(t-2\epsilon)  \Tr_\G \left( tD_0 \dd{u}(\tD_u)_{|u=0} \tD_0
e^{-t\tD_0^2} - \tD_0^2 \dd{u}(\tD_u)_{|u=0} 
e^{-t\tD_0^2} \right)\nonumber\\
& & {}  = -(t-2\epsilon)  \Tr_\G \left(\dd{u}(\tD_u)_{|u=0} \tD_0^2 
e^{-t\tD_0^2} \right) \textrm{ by permuting factors as in (\ref{2})}.
\end{eqnarray*}
Taking the limit $\epsilon \rightarrow 0$, we get
\begin{eqnarray}
\label{5}
-2t \Tr_\G \left(\dd{u}(\tD_u)_{|u=0} \tD_0^2 e^{-t\tD_0^2} \right) 
  = 2t \dd{t}\Tr_\G \left(\dd{u}(\tD_u)_{|u=0}e^{-t\tD_0^2} \right). 
\end{eqnarray}
To make the corresponding evaluation for the left-hand side of (\ref{3}) note
that  
\begin{eqnarray}
\label{6}
\dd{u} \left( \Tr_\G \left( \tD_u e^{-t\tD_u^2}\right ) \right) = \Tr_\G \left( \dd{u}(\tD_u)e^{-t\tD_u^2}\right ) + \Tr_\G \left( \tD_u \dd{u}(e^{-t\tD_u^2}) \right ).
\end{eqnarray}
Also for fixed $\epsilon$, 
as $e^{-t\tD_u^2}= 
e^{-(t-\epsilon)\tD_u^2}e^{-\epsilon \tD_u^2}=
e^{-\epsilon \tD_u^2}e^{-(t-\epsilon)\tD_u^2}$,
{\setlength\arraycolsep{2pt}
\begin{eqnarray}
\label{7}
\dd{u} \Tr_\G \left( \tD_u e^{-t \tD_u^2} \right) &=&  
\Tr_\G \left( \dd{u} (\tD_u) e^{-t \tD_u^2} \right)  \nonumber\\
& & {} + \Tr_\G \left(\tD_u \dd{u}(e^{-(t-\epsilon)\tD_u^2})
e^{-\epsilon \tD_u^2} \right) \nonumber\\
& & {}  + \Tr_\G \left(\tD_u e^{-(t-\epsilon)\tD_u^2}
\dd{u}(e^{-\epsilon \tD_u^2}) \right).  
\end{eqnarray}}
Letting $\epsilon \rightarrow 0$ and comparing (\ref{6}) with (\ref{7}) gives
\begin{eqnarray*}
\label{8}
\lim_{\epsilon \rightarrow 0} \Tr_\G \left(\tD_u e^{-(t-\epsilon)\tD_u^2}
\dd{u}(e^{-\epsilon \tD_u^2}) \right) = \lim_{\epsilon \rightarrow 0} \Tr_\G \left(\tD_u \dd{u}(e^{-\epsilon \tD_u^2})e^{-(t-\epsilon)\tD_u^2} \right)=0. 
\end{eqnarray*}
So the derivative at $u=0$ of the left-hand side of (\ref{3}) is
\begin{eqnarray}
\label{9}
\lim_{\epsilon \rightarrow 0} \left[ \Tr_\G \left( \tD_0 \dd{u} (e^{-(t-\epsilon) \tD_u^2})_{|u=0} e^{-\epsilon \tD_0^2}
\right) - \Tr_\G \left(  \tD_0 \dd{u}(e^{-\epsilon \tD_u^2})_{|u=0} e^{-(t-\epsilon) \tD_0^2}\right) \right]  \nonumber\\
{} = \dd{u} \Tr_\G \left( \tD_u e^{-t\tD_u^2}\right)_{|u=0} - \Tr_\G \left( \dd{u} (\tD_u)_{|u=0} e^{-t\tD^2_0} \right).
\end{eqnarray}
Combining (\ref{5}) and (\ref{9}) gives
\begin{eqnarray*}
\label{10}
\dd{u} \left( \Tr_\G \left( \tD_u e^{-t \tD_u^2} \right)\right)_{|u=0} =
\Tr_\G \left( \dd{u} (\tD_u)_{|u=0} e^{-t\tD_0^2} \right) 
+ 2t \dd{t} \Tr_\G \left( \dd{u}(\tD_u)_{|u=0} e^{-t \tD_0^2} \right). 
\end{eqnarray*}

This ends the proof of the formula (\ref{Formula}).

If $\dd{u} \eta_\G(u)$ exists, using the equality (\ref{Formula}),
it would be equal to the following limits
{\setlength\arraycolsep{2pt}
\begin{eqnarray}
\label{Step-a}
&& {}\lim_{T \rightarrow \infty}  
\lim_{\epsilon \rightarrow 0}
\left( \frac{1}{\G(1/2)} \int_\epsilon^T \!\! t^{-1/2} \Tr_\G \left( 
(\dd{u} \tD_u) e^{-t\tD_u^2} \right) \d t \right.
\nonumber\\
& & {} +\left. \frac{1}{\G(1/2)} \int_\epsilon^T \!\! 2  t^{1/2} 
\dd{t} \left( \Tr_\G
\left( (\dd{u} \tD_u) e^{t \tD_u^2} \right) \right) \d t
 \right).
\end{eqnarray}}
Now, in the previous expression (\ref{Step-a}), 
we integrate the second term by parts 
and  obtain
\begin{equation*}
\frac{1}{\G(1/2)} \lim_{T \rightarrow \infty}
 2 T^{1/2} \Tr_\G \left(
\dd{u} (\tD_u) e^{-T \tD_u^2} \right)  - \frac{1}{\G(1/2)}
 \lim_{\epsilon \rightarrow 0} 2 e^{1/2} \Tr_\G \left(
\dd{u} (\tD_u) e^{-\epsilon \tD_u^2} \right). 
\end{equation*}
We will now compute the limits of these both terms.

($b$) We will show that:
\begin{eqnarray}
\label{Zero}
\lim_{T \rightarrow \infty} T^{1/2} \Tr_\G \left(
\dd{u} (\tD_u) e^{-T \tD_u^2} \right) = 0.
\end{eqnarray}

For this, we recall the useful Proposition $4.15$ of \cite{Ati76}:
for a bounded, $\G$-invariant operator $A$ on $L^2(\tilde{S})$ 
which is of $\G$-trace class and  a sequence of bounded, 
$\G$-invariant operators  $B_j$ on $L^2(\tilde{S})$ which converges to a 
bounded, $\G$-invariant operator $B$ on $L^2(\tilde{S})$ i.e.
$\forall f \in L^2(\tilde{S}),\;\; B_j f \rightarrow B f$ in $L^2$-norm, we
have $ \displaystyle \lim_{ j\rightarrow + \infty} \Tr_{\G}(AB_j) =
\Tr_\G(AB)$.  

Let $T_0$ be a strictly positive real.
Consider $$A = \dd{u} (\tD_u) e^{-T_0 \tD_u^2} = \left(\dd{u} (\tD_u) 
e^{-\frac{T_0}{2} \tD_u^2}\right)
\left(e^{-\frac{T_0}{2} \tD_u^2}\right).$$ 
Since the operator $A$ is the composition
of two $\G$-Hilbert-Schmidt operators (Sub-lemma \ref{bounded}), 
the operator $A$ is, by
definition, of $\G$-trace class. 
Let $T_j$ be an increasing sequence
of real which tends to infinity. We defined
$$B_j=  T_j^{1/2}e^{-(T_j-T_0) \tD_u^2}= T_j^{1/2}e^{-(T_j-T_0)c^2} 
e^{-(T_j-T_0) (\tD_u^2-c^2)}$$
 where $c$ is chosen such that $\tD_u^2-c^2 > 
\alpha >0$ (it's possible because $0\notin \spec (\tD_u)$).
Then the spectral
theory shows that  $e^{-(T_j-T_0) (\tD_u^2-c^2)}$ converges to $B$ the operator
of projection on the kernel of $\tD_u^2-c^2$ which is $\{0\}$. On the
other hand $T_j^{1/2}e^{-(T_j-T_0)c^2}$ converges to $0$ as $j$ tends to
infinity. Therefore we have proved that  $B_j$ converges to the 
bounded $\G$-invariant operator $B$ (projection on $\{0\}$).
Finally we have the desired formula (\ref{Zero}).

($c$) Now we prove the formula:
\begin{eqnarray}
\label{Local}
 - \frac{1}{\G(1/2)} \lim_{\epsilon \rightarrow 0} 2 \epsilon^{1/2} \Tr_\G 
\left( \dd{u} (\tD_u) e^{-\epsilon \tD_u^2} \right) = \int_F \!\! 
a_n(x, \dd{u} \tD_u, \tD_u) \d x
\end{eqnarray}
where $F$ is a fundamental domain and $a_n$ is a local invariant in the jets
of the symbols of $(\dd{u} \tD_u,\tD_u)$.

Using  symbolic calculus 
(see  \cite[$\S 1.7$]{Gil84}), 
 we can improve Sub-lemma \ref{bounded} as follows.

\begin{sublemma}
Under the same assumptions as in  Sub-lemma \ref{bounded},
there is an asymptotic expansion on the diagonal:
$$k(t,\tilde{x},\tilde{x})=
\{\tilde{Q}_xK(t,\tilde{x},\tilde{y})\}_{\tilde{y}=\tilde{x}} \sim 
\sum_{k=0}^{\infty} t^{(k-n-a)/2}e_k(\tilde{x},\tilde{Q},\tD)$$
where $e_k$ are smooth local invariants of the jets of the symbols of $\tD$ 
and $\tilde{Q}$.
This means that given any integer $l$ there exists an integer $k(l)$ such that
$$|k(t,\tilde{x},\tilde{x})- \sum_{0 \leq k \leq k(l)} t^{(k-n-a)/2}e_k(\tilde{x})| <
C_kt^k \textrm{ for } 0<t<1.$$
\end{sublemma} 

\begin{SketchProof}
The proof of this sub-lemma is nearly the same as when the manifold
is closed, which is the case in  \cite[Lemma $1.7.7.$]{Gil84}.
As in the proof of Sub-lemma \ref{bounded}, we use the non-compact
Sobolev embedding associated to the elliptic estimates in order to
compute the norm of the operators from $(L^2(\tilde S), \|.\|)$
to $(\CB^r(\tilde S), \|.\|_{\CB^r})$. Here the fact that we are on a 
covering of a closed manifold with lifted operators plays
an important role. In the closed manifold case, we use the
norm of the operators from the Sobolev space $(L^2(S), \|.\|)$
onto higher order Sobolev spaces and use
the Sobolev embedding only at the end of the proof. 
The remainder of the proof 
is local 
and thus is exactly the same.
\end{SketchProof}

If we have such an asymptotic expansion for the kernel, we deduce immediately
the following asymptotic expansion for the $\G$-trace:
$$ \Tr_\G \left(\tilde{Q} e^{-t \tD} \right) =
\int_F \Tr k(t,\tilde{x}, \tilde{x}) \, \d \vol (\tilde{x}) \sim 
\sum_{k=0}^{\infty} t^{(k-n-a)/2} \int_F \Tr e_k(\tilde{x},\tilde{Q},\tD) \,
\d \vol (\tilde{x}).$$
And therefore, in our case, with $\tilde{Q}=\dd{u}(\tD_u)$ and 
$\tilde{D}=\tD_u$, 
we obtain
$$ \epsilon^{1/2} \Tr_\G \left(\dd{u}(\tD_u) e^{-\epsilon \tD_u^2} \right) 
 \sim 
\sum_{k=0}^{\infty} \epsilon^{(k-n)/2} \int_F \Tr e_k(\tilde{x},\dd{u}(\tD_u),
\tD_u) \,
\d \vol (\tilde{x}).$$
As the symbols of the operators and of their lifts are the same, we can replace
the integration over the fundamental domain $F$ by an integration over $V$,
thus
$$ \epsilon^{1/2} \Tr_\G \left(\dd{u}(\tD_u) e^{-\epsilon \tD_u^2} \right) 
 \sim 
\sum_{k=0}^{\infty} \epsilon^{(k-n)/2} \int_V \Tr e_k(x,\dd{u}(\D_u),
\D_u) \,
\d \vol (x).$$
We know that $\dd u \bar \eta $ exists (see \cite[$\S 1.10$]{Gil84}).
By the same computation as for $\dd u \bar \eta_\G $, we show that
the limit 
$\lim_{\epsilon \rightarrow 0} \epsilon^{1/2} 
\Tr_\G \left(\dd{u}(\D_u) e^{-\epsilon \D_u^2} \right)$
is finite and is equal to $-\frac{1}{2} \Gamma(1/2) \dd u \bar \eta $.
Moreover, we have also
$$ \epsilon^{1/2} \Tr_\G \left(\dd{u}(\D_u) e^{-\epsilon \D_u^2} \right)
 \sim 
\sum_{k=0}^{\infty} \epsilon^{(k-n)/2} \int_V \Tr e_k(x,\dd{u}(\D_u),
\D_u) \,
\d \vol (x).$$
As a consequence 
we obtain  $\int_F \Tr e_k(\tilde{x},\dd{u}(\tD_u) =
\int_V \Tr e_k(x,\dd{u}(\D_u) 
=0 $ for $k<n$.
The only terms which arise in the computation of the limit (\ref{Local}), 
are the ones for $k \geq n$, which ones are
non-zero when $\epsilon \rightarrow 0$ only 
for $k=n$ and the result  follows.

 Finally we combine the points ($a$),($b$) and ($c$) to achieve 
the proof of Lemma \ref{EtaLocal}. 
\end{proof}

\section{Approximate inequality}
\label{RKato}

In the present section we will develop the second step of the proof 
as announced in  Section \ref{Outline}.

\subsection{Refined Kato inequality}
First we introduce the context.
Let $V$ be a closed Riemannian spin manifold of dimension $n$.
Let $\D: C^{\infty}(S) \rightarrow C^{\infty}(S)$ be a first order operator.
The \emph{twisted operator}
$\D^Y$ of $\D$ by a unitary bundle $(Y,\nabla^Y)$ is the operator whose 
symbol is $\sigma(\D^Y)=\sigma(\D) \otimes 1$ defined 
as follows
$$
\begin{array}{ccccc}
C^{\infty}(S\otimes Y) & \stackrel{\nabla^{S \otimes Y}}{\longrightarrow} & 
C^{\infty}(T^*V \otimes (S \otimes Y)) =
C^{\infty}((T^*V \otimes S) \otimes Y) & 
\stackrel{ \sigma(\D)\otimes 1}{\longrightarrow} & C^{\infty}(S\otimes Y)\\
\end{array}$$
So we have $\D^Y= (\sigma(\D) \otimes 1) \circ \nabla^{S \otimes Y} =
\sigma \circ \nabla$ where $\sigma$ is just the
symbol of the operator.  
We now recall the refined Kato inequality \cite{CGH00}: 

\begin{proposition}
\label{Kato}
Consider $\Phi$ an element of $\ker \sigma$ at some point and $\phi$ 
an element of 
$S \otimes Y$ at the same point.
Then
$$ \sup_{|\phi|=1} |\langle \Phi, \phi \rangle | \leq k_\sigma |\Phi|,$$
where the constant $k_\sigma (<1)$ only 
depends on the symbol of the operator $\D$. \\
Futhermore,
for $\psi \in \ker \D^Y$ i.e. $\Phi = \nabla \psi \in \ker \sigma$, 
$|\d |\psi|| \leq k_\sigma |\nabla \psi|$ $(*)$.

In particular,
for any (twisted) Dirac operator $\Dirac^Y$, 
we know $k_\sigma=
((n-1)/n)^{1/2}$. Moreover
we have equality in the previous inequality $(*)$ if and only if
there exists 
a $1$-form $\alpha$ such that $$\nabla \psi = \alpha \otimes \psi +  
\frac{1}{n} \sum e_i \otimes e_i \cdot \alpha
\cdot \psi.$$
\end{proposition}

\begin{proof} 
We adapt the general proof to the particular case of a Dirac operator.
It is well known that the symbol $\sigma$ of the Dirac operator is 
the Clifford multiplication.
We define 
$\Pi : C^{\infty}(T^*V \otimes (S \otimes Y)) \rightarrow 
C^{\infty}(T^*V \otimes (S \otimes Y))$ by
$\Pi (\alpha \otimes \psi) = -\frac{1}{n} \sum e_i \otimes e_i \cdot \alpha
\cdot \psi$ where $e_i$ is an orthonormal basis of $T^*V$, $\alpha \in T^*V$,
$\psi \in S \otimes Y$ and $\cdot$ the Clifford multiplication. 
We have $\ker \sigma = \ker \Pi$, the endomorphisms $\Pi$ and $1-\Pi$ are
orthogonal projections. Let $\phi$ be a section of $S \otimes Y$ 
with norm $1$. Thus, if $\Phi \in \ker \Pi$ we have
$|\langle \Phi, \phi \rangle | = \sup_{|\alpha|=1} |\langle \Phi, 
\alpha \otimes \phi  \rangle |$. But, by the Cauchy-Schwarz inequality,
$|\langle \Phi, 
\alpha \otimes \phi  \rangle |= |\langle \Phi, 
(1-\Pi)(\alpha \otimes \phi)  \rangle | \leq |\Phi| |(1-\Pi)
(\alpha \otimes \phi)| $. 
We now write Pythagoras' equality
 $1=|\alpha \otimes \phi|^2 = |\Pi(\alpha \otimes \phi)|^2 +|(1-\Pi)
(\alpha \otimes \phi)|^2 $.  
We can compute directly 
\begin{equation*}
|\Pi(\alpha \otimes \phi)|^2  =  \frac{1}{n^2} \langle
\sum_i e_i \otimes e_i \cdot \alpha \cdot \phi,  
\sum_j e_j \otimes e_j \cdot \alpha \cdot \phi \rangle= \frac{1}{n} |\alpha|^2 |\phi|^2 = \frac{1}{n}.
\end{equation*}
The constant $k_\sigma$ follows.

As in \cite[Theo. $3.1$]{CGH00}, we can study the equality case.
In this case, we have equality in the Cauchy-Schwarz inequality, hence
there exists a $1$-form $\alpha$ and a section $\phi$ of $S \otimes Y$
such that  $\Phi=(1-\Pi)(\alpha \otimes \phi)$. Furthermore, for 
$\psi \in \ker \Dirac^Y$, the equality implies that there exists 
a $1$-form $\alpha$ such that $\nabla \psi = \alpha \otimes \psi +  
\frac{1}{n} \sum e_i \otimes e_i \cdot \alpha
\cdot \psi$. This will be useful in Section \ref{Proof1}.
\end{proof}

\subsection{Bochner-Lichnerowicz-Weitzenb\"ock formula}
For this section we refer to the article \cite[\S $2$]{GL83}.
Let $\D$ denote a Dirac operator acting on a Dirac bundle $S$ on 
$V$ and $\tD$
denote the lifted Dirac operator acting on $\L^2(\tilde{S})$, 
the Hilbert space of $\L^2$-integrable sections of $\tilde{S}$ on $\tV$.
The inner product and the norm on $\L^2(\tilde{S})$ will be denoted by
$\lps .,.\rps$, respectively $\|.\|$ ; they are induced as usual by the 
pointwise inner product $\langle .,. \rangle$ on $\tilde{S}$ and
the canonical volume form $\d \vol_{\tV}$ on $\tV$.

We define the \emph{Sobolev space} $\L^{1,2}(\tilde{S})$. 
For
$\psi \in C^\infty_c(\tilde{S})$, the smooth sections of $\tilde{S}$ 
which have compact support, we set
$$\|\psi\|^2_1= \int_{\tV} (\langle \psi,\psi \rangle + \langle \tnabla\psi, \tnabla\psi
\rangle)$$
and denote $\L^{1,2}(\tilde{S})$ the completion of
$C^\infty_c(\tilde{S})$ in this norm.

Let consider 
the twisted Dirac operators. Let $Y$ be an hermitian bundle 
over $V$ and $\Dirac^Y$ be 
the associated twisted 
Dirac on $V$ acting on sections of $\Spin \otimes Y$.
Let $\tilde{\Dirac^Y}$ be the lifted Dirac operator on $\tV$ and 
$\tilde{Y}$ be the lifted bundle over $\tV$. The operator  $\tilde{\Dirac^Y}$
is the same as the Dirac operator $\tDirac$ twisted by the bundle $\tilde{Y}$.
 We emphasize that,
on $\tV$ (as on $V$), we have the fundamental 
Bochner-Lichnerowicz-Weitzenb\"ock formula 
$$\|\tilde{\Dirac^Y} \psi \|^2 = \| \tilde{\nabla} \psi \|^2 + \lps 
\frac{\scal}{4} \psi,\psi \rps + \lps \R^{\tilde{Y}} \psi, 
\psi \rps $$
where $\R^{\tilde{Y}}$ is defined by the formula
$$\R^{\tilde{Y}}(s \otimes e) =
\frac{1}{2} \sum_{i,j} (e_i \cdot e_j \cdot s) \otimes
 R^{\tilde{Y}}_{e_i,e_j}(e) \textrm{ where } R^{\tilde{Y}}_{v,w} = 
[\tnabla^{\tilde{Y}}_v, \tnabla^{\tilde{Y}}_w]-\tnabla^{\tilde{Y}}_{[v,w]}$$
with $s$ a section of $\tilde{\Spin}$, $e$ a
section of $\tilde{Y}$, $(e_i)$ an orthonormal basis of $T^*\tV$,
$\cdot$ is the Clifford multiplication
and $R^{\tilde{Y}}$ is the curvature $2$-form with value
in $\operatorname{End}(\tilde{Y})$. Moreover, we can bound above the
norm of the operator $\R^{\tilde{Y}}$ as follows:
\begin{equation}
\label{CurvatureTerm}
\| \R^{\tilde{Y}} \| \leq \alpha_{n} \| R^{\tilde{Y}} \|
= \alpha_{n} \| R^{Y} \| 
\end{equation}
where $\alpha_{n}$ only depends on the dimension 
of $V$.

As a consequence of the  Bochner-Lichnerowicz-Weitzenb\"ock formula 
and of the fact that  the curvature term 
$\scal/4 + \R^{\tilde{Y}}$ is uniformly bounded on $\tV$, the
maximal domain of $\tDirac$ on $\L^2(\tilde{\Spin})$ is exactly 
$\L^{1,2}(\tilde{\Spin})$ (see \cite[Theo. 2.8]{GL83}).

\subsection{Key Lemma}

The main purpose of the section is to prove the following Lemma.

\begin{lemma}
\label{KeyLemma}
\label{UKato}
Let $Y$ be a Riemannian bundle over $V$ and 
$\Dirac^Y$ be the associated twisted 
Dirac operator on $V$. Let $\tilde{\Dirac^Y}$ 
be the lifted Dirac operator on $\tV$.
If $0$ belongs to the spectrum of $\tilde{\Dirac^Y}$, we have
$$\inf \scal(V,g) \leq -4 \frac{n}{n-1} \lambda_0(\tV, \tilde{g})
+\alpha_n \|R^{Y}\|$$
where $\alpha_n$ only depends on the dimension of $V$.
\end{lemma}

\begin{proof}
We need to distinguish two cases. 
The trivial one is when there exists an eigenvector associate to the 
eigenvalue $0$. Then we  write the Bochner formula and can conclude 
in a trivial way. If not, $0$ is in the essential spectrum and we 
cannot directly use the refined Kato inequality.

For sake of simplicity
we will denote $\D$ for $\Dirac^Y$ and $S$ for $\Spin$.

\emph{If $0$ is in the point spectrum} there exists an eigenspinor $\psi \in 
\L^{1,2}(\tilde{S}^+ \otimes \tilde{Y})$ of the operator $\tD$ on $\tV$. 
We apply the Bochner-Weitzenb\"ock formula to this $\psi$:
$$0=\|\tD \psi \|^2 = \| \tilde{\nabla} \psi \|^2 + \lps 
\frac{\scal}{4} \psi,\psi \rps + \lps \R^{\tilde{Y}} \psi, 
\psi \rps. $$
To obtain the desired inequality, we estimate the three terms of this equality.
For the first term, we use successively
the refined Kato inequality  (Proposition \ref{Kato}) and the  
Rayleigh's quotient. For the second,
 we just bound below the scalar curvature with  its infimum. 
For the third, we use the upper bound (\ref{CurvatureTerm}) of the norm of
the operator $\R^{\tilde{Y}}$.

We finally obtain
$$ 0 \geq \left(\frac{n}{n-1} \lambda_0 (\tV,\tilde{g}) + 
\frac{\inf \scal(V,g)}{4} -\alpha_n \|R^Y\|\right) 
\|\psi\|^2.$$ 
As $\|\psi\|^2 >0$ the result follows.

\emph{If $0$ is in the essential spectrum}, there exists a sequence $\psi_k \in
\L^{1,2}(\tilde{S}^+ \otimes \tilde{Y})$ such that 
$\|\psi_k\|^2=1$ and $\| \tD \psi_k\|^2 \leq \frac{1}{k}
\|\psi_k\|^2 = \frac{1}{k}$.
We have $2 |\d |\psi|| \, |\psi| = |\d |\psi|^2| = 2 |\langle \tnabla \psi, \psi 
\rangle |$, thus $| \d |\psi||^2 = | \langle \tnabla \psi, \phi \rangle|$ 
where $\phi= \frac{\psi}{|\psi|}$.
Writing   $\tD = \Pi \circ \tnabla$, we obtain:
$$|\langle \tnabla \psi, \phi \rangle |^2 = |\langle \Pi \circ \tnabla \psi, 
\phi
 \rangle|^2 +|\langle (1-\Pi) \circ  \tnabla \psi, \phi \rangle|^2.$$
As $\Phi = (1-\Pi) \circ \tnabla \psi \in \ker \Pi$, $|\langle \Phi, \phi 
\rangle| \leq k_\Pi |\Phi| \leq k_\Pi |\tnabla \psi|$ (Proposition \ref{Kato}).
 Moreover 
$|\langle \Pi \circ \tnabla \psi, \phi \rangle| = 
|\langle \tD \psi, \phi \rangle| 
\leq |\tD \psi||\phi|=|\tD \psi|$, hence
$$|\langle \tnabla \psi, \frac{\psi}{|\psi|} \rangle|^2 \leq |\tD \psi|^2 + 
k_\Pi^2 |\tnabla \psi |^2.$$
Applying this to $\psi_k$ and replacing $k_\Pi$ by its value 
$\sqrt{\frac{n-1}{n}}$,
 we obtain
$$0 \leq \int |\d |\psi_k||^2 \leq \frac{1}{k} + \frac{n-1}{n}
\|\tnabla \psi_k \|^2.$$

Now we write the Bochner-Weitzenb\"ock formula applied to $\psi_k$
$$\| \tD \psi_k\|^2 = \|\tnabla \psi_k\|^2 + 
\lps \frac{1}{4} \scal \psi_k ,\psi_k \rps
+\lps \R^{\tilde{Y}} \psi_k, \psi_k \rps,$$
consequently
$$\frac{1}{k} \geq \frac{\int |\d |\psi_k||^2 - \frac{1}{k}}
{\frac{n-1}{n}} + \frac{1}{4} \inf \scal \|\psi_k \|^2 
- \alpha_n \|R^Y\| \|\psi_k\|^2.$$
Now we use the Rayleigh's quotient
$$\frac{1}{k} \geq \frac{\lambda_0(\tV) - \frac{1}{k}}
{\frac{n-1}{n}} + \frac{1}{4} \inf \scal  
- \alpha_n \|R^Y\| .$$
Letting $k \rightarrow \infty$
$$0 \geq \frac{\lambda_0(\tV)}
{\frac{n-1}{n}} + \frac{1}{4} \inf \scal  
- \alpha_n \|R^Y\| .$$
This ends the proof of Lemma \ref{UKato}.
\end{proof}

\section{Proof of Theorem \ref{AgenusCase} and Theorem \ref{EqualityCase}}
\label{Proof1}

\subsection{Proof of Theorem \ref{AgenusCase}}

Let $V$ be a closed Riemannian spin manifold with non-zero 
$\Ah$-genus.
We use the Atiyah $\G$-index Theorem  
to ensure that
$\ind_\G \tDirac^+ = \ind \Dirac^+ = \Ah_V[V] \neq 0$ thus there exists $\psi 
\in \L^{1,2}(\tilde{S}^+), \; \psi \neq 0$ such that $\tDirac \psi = 0$.
Therefore we can use the key Lemma \ref{UKato} 
simply for the (untwisted) usual Dirac operator
and obtain directly the desired inequality.

\subsection{Equality case (Theorem \ref{EqualityCase})}

As we know the equality case in the refined Kato inequality 
\cite[Theo. $3.1$]{CGH00}, we 
can study the equality case in Theorem \ref{AgenusCase}.
Since one of the implications is trivial (i.e. $\scal (V,g) = 0$ implies
the equality), we have only to prove the other one. So we state the following
Proposition. 

\begin{proposition}
\label{FirstStep}
Under the same assumptions as in Theorem \ref{AgenusCase}. If we are 
in the equality case, that is
$$\inf \scal(V,g) = -4 \frac{n}{n-1} \lambda_0(\tV, \tilde{g}),$$
then  $\tV$ is compact and $\scal (V,g) = 0$.
\end{proposition}

\begin{proof}
Suppose that we are in the equality case
of Theorem \ref{AgenusCase}, we distinct two cases: $\tV$ is compact or
$\tV$ is non compact. If $\tV$ is compact we are going to
show that $\scal(V,g)=0$. On the
other hand, we will show that the case $\tV$ non-compact is impossible.

\emph{If $\tV$ is compact,} then we have trivially 
$\lambda_0(\tV ,\tilde{g})=0$ and hence we get 
$\inf \scal(V,g) = -4 \frac{n}{n-1} 
\lambda_0(\tV ,\tilde{g})=0$. Using the $\widehat{A}$-genus vanishing Theorem
of Lichnerowicz we conclude that $\scal(V,g)=0$.

\emph{When $\tV$ is non-compact,} we will find a contradiction.

We begin with the equality cases in  
the inequalities that we use   
in the proof of the 
key Lemma \ref{UKato}:
\begin{itemize}
\item[($a$)] Equality in the Refined Kato inequality: 
$ | \d |\psi||^2 = \frac{n}{n-1} |\tnabla \psi |^2$
and, as recalled in Proposition \ref{Kato}, 
there exists a $1$-form $\alpha$ on $\tV$ such that  
$  \tnabla \psi = \alpha \otimes \psi + \frac{1}{n} \sum_i e_i 
\otimes e_i \cdot \alpha \cdot \psi$, 
where $e_i$ is an orthonormal basis of $T^*\tV$ and
$\cdot$ is the Clifford multiplication,

\item[($b$)] Equality for the Rayleigh's quotient: 
 $\lambda_0(\tV) = \frac{ \int | \d |\psi||^2}{\int |\psi|^2}$ 
therefore $\Delta |\psi| = \lambda_0(\tV) |\psi|$, 

\item[($c$)] Equality for the infimum: 
 $ \inf \scal(V,g) = \scal(V,g)$ and thus
$\lambda_0 (\tV,\tilde{g}) = -\frac{n-1}{n} \frac{\scal(V,g)}{4}$.
\end{itemize}

Next we make some intermediate computations in order to conclude to a 
contradiction in a last time.

With ($a$) we compute $|\tnabla \psi|^2$ as follows
\begin{eqnarray*}
|\tnabla \psi |^2 &=& \langle \alpha \otimes \psi + \frac{1}{n} \sum_i e_i 
\otimes e_i \cdot \alpha \cdot \psi,
\alpha \otimes \psi + \frac{1}{n} \sum_j e_j 
\otimes e_j \cdot \alpha \cdot \psi \rangle \nonumber\\
&=& |\alpha|^2|\psi|^2 + \frac{2}{n} \sum_i  \langle e_i,\alpha\rangle
 \Re \langle \psi , e_i \cdot \alpha \cdot \psi \rangle + \frac{1}{n^2}
\sum_{i,j} \langle e_i, e_j \rangle \langle e_i \cdot \alpha \cdot \psi,
 e_j \cdot \alpha \cdot \psi\rangle \nonumber \\
&=& |\alpha|^2|\psi|^2 + \frac{2}{n} \sum_{i,j}  \langle e_i,\alpha\rangle
\langle e_j,\alpha\rangle
\Re \langle \psi , e_i \cdot e_j \cdot \psi \rangle +  \frac{1}{n^2}
\sum_i \langle e_i \cdot \alpha \cdot \psi,
 e_i \cdot \alpha \cdot \psi\rangle \nonumber\\
&=& |\alpha|^2|\psi|^2 + \frac{2}{n} \sum_i \langle e_i,\alpha\rangle^2
\Re \langle \psi, e_i \cdot e_i \cdot \psi \rangle \nonumber\\
 &&+\frac{2}{n}
 \sum_{i<j} \langle e_i,\alpha\rangle  \langle e_j,\alpha\rangle
 \Re \langle \psi , (e_i \cdot e_j+ e_j \cdot e_i) \cdot \psi \rangle
  +  \frac{1}{n^2}
\sum_i |\alpha|^2 |\psi|^2. 
\end{eqnarray*}
Hence, we obtain
\begin{eqnarray}
|\tnabla \psi |^2 &=& \frac{n-1}{n}|\alpha|^2|\psi|^2.
\label{a}
\end{eqnarray}

We will now  find the contradiction proving first:

\begin{lemma}
\label{Delta}
We have $|\alpha|^2=  \left(\frac{n}{n-1}\right)^2 \lambda_0(\tV, \tg)$ 
and $\Delta |\psi|^2=0$.
\end{lemma}

\begin{proof}
We have  two different manners to calculate
 $\frac{1}{2} \Delta |\psi|^2$:

\begin{itemize}
\item[$\bullet$] First,
\begin{eqnarray*}
 \frac{1}{2} \Delta |\psi|^2 & = & - |\tnabla \psi|^2 - 
\frac{1}{4} \scal |\psi|^2  \\
&\stackrel{(\ref{a})}{=}& -\frac{n-1}{n} |\alpha|^2 |\psi|^2 -
\frac{1}{4} \scal |\psi|^2 \\ 
&\stackrel{(c)}{=}& \left(- \frac{n-1}{n} |\alpha|^2 + \frac{n}{n-1}
\lambda_0 \right) |\psi|  
\end{eqnarray*}
where we have used
$ \frac{1}{2} \Delta |\psi|^2 = - |\tnabla \psi |^2-\frac{1}{4} 
\scal |\psi|^2.$ Indeed
we always have $ \frac{1}{2} \Delta |\psi|^2  = 
\langle \tnabla^* \tnabla \psi, \psi \rangle - |\tnabla \psi|^2$.
As $\psi$ is harmonic, $0=\langle \Dirac^2\psi, \psi \rangle = \langle 
\tnabla^* \tnabla 
\psi + \frac{1}{4} \scal \psi, \psi \rangle$, hence we have 
the desired formula.

\item[$\bullet$] On the other hand, since
$ \Delta f^2 = - 2 |\d f|^2 + 2f \Delta f,$ we have
\begin{eqnarray*} 
 \frac{1}{2} \Delta |\psi|^2 &=& -|\d |\psi||^2 + |\psi|
 \Delta |\psi|  \\
& \stackrel{(a)}{=} & -\frac{n-1}{n} |\tnabla \psi|^2+|\psi| \Delta |\psi|  \\ 
& \stackrel{(\ref{a})}{=} &  -\left(\frac{n-1}{n}\right)^2 |\alpha|^2
 |\psi|^2 +|\psi| \Delta  |\psi|   \\
& \stackrel{(b)}{=} & -\left(\frac{n-1}{n}\right)^2 |\alpha|^2
 |\psi|^2 + \lambda_0 |\psi|^2  \\ 
& = & \left(-\left(\frac{n-1}{n}\right)^2 |\alpha|^2  + \lambda_0\right)
|\psi|^2.  
\end{eqnarray*}
\end{itemize}
Finally we have
$$-\left(\frac{n-1}{n}\right)^2 |\alpha|^2  + \lambda_0 =- \frac{n-1}{n} |\alpha|^2 + \frac{n}{n-1} \lambda_0$$
consequently 
$$|\alpha|^2 = \left(\frac{n}{n-1}\right)^2 \lambda_0.$$
Thanks to this equality, we deduce that $\Delta |\psi|^2 =0$.
\end{proof}

On the other hand we check that
\begin{equation}
\label{Deriv}
\d |\psi|^2 = 2 \frac{n-1}{n} \alpha |\psi|^2. 
\end{equation}
Indeed we have
$$\d |\psi|^2= 2 \Re \langle \tnabla \psi, \psi \rangle \stackrel{(a)}{=}
 2 \alpha
|\psi|^2+ \frac{2}{n} \sum_i \Re \langle e_i \cdot \alpha \cdot \psi,\psi 
\rangle e_i.$$
For $X \in \ker \alpha$, we have 
$\sum_i \Re \langle e_i \cdot \alpha \cdot \psi,\psi 
\rangle e_i(X) = \Re  \langle X^\flat \cdot \alpha \cdot \psi, \psi\rangle =0$
using the properties of the Clifford algebra and the fact that
the  hermitian
metric is adapted to the Clifford multiplication.
For $X = \alpha ^\#$, we compute 
$\sum_i \Re \langle e_i \cdot \alpha \cdot \psi,\psi 
\rangle e_i(\alpha ^\#) = \Re  \langle \alpha \cdot \alpha \cdot \psi, 
\psi\rangle = - |\alpha|^2|\psi|^2$.
The equality (\ref{Deriv}) follows.

We now  improve Lemma \ref{Delta} as follows:
\begin{lemma}
\label{AlphaZero}
We have $|\alpha|= 0$ and so $\lambda_0(\tV, \tg) =0$.
\end{lemma}

\begin{proof}
\emph{By contradiction}, we suppose that the constant $|\alpha|$ is non-zero 
i.e. $\lambda_0 \neq 0$. Let $f=|\psi|^2$. We know that $f \in \L^1(\tV)$
(because $\psi \in \L^2(\tilde S)$),
$\Delta f =0$ (Lemma \ref{Delta}) and $|\grad f|=2 \sqrt{\lambda_0}f$ (Equality
(\ref{Deriv}) and Lemma \ref{Delta}). Notice that
$f > 0$ everywhere because $\psi$ is a non-zero harmonic spinor.  

We define the vector field 
$$X= \frac{\grad f}{|\grad f|}$$ 
which is the unit normal to the level-submanifolds $\tV_c=f^{-1}(c)$. 
The vector field $X$ is complete 
on $\tV$ because $|X|$ is constant, equal to $1$.

On an integral curve $\gamma$ of $X$, $\dd{t} f =2 \sqrt{\lambda_0}f$, 
therefore $f(t) = k_\gamma e^{2 \sqrt{\lambda_0} t}$ where $k_\gamma$ is a 
strictly positive constant depending
on the curve.

If there exists a closed integral curve $\gamma$ 
then there are two  different times 
$t$ and $t'$  such
that $k_\gamma e^{2 \sqrt{\lambda_0}t}=k_\gamma e^{2 \sqrt{\lambda_0}t'}$: 
it is impossible. Therefore
no such curve is  closed. 

Choose $f_0 \in \Im f$ and 
a little transverse neighbourhood $\EuScript{V}_0 \subset \tV_{f_0}$ of the
vector field $X$ 
around a point $x_0 \in \tV_{f_0}$. Let $\gamma$ the integral curves of
$X$ such that $\gamma(0) \in \EuScript{V}_0$ thus 
we have $f(t)=f_0 e^{2 \sqrt{\lambda_0}t}$ on this curves.
 The  neighbourhood $\EuScript{V}_0$ is 
chosen sufficiently small such that the flow of $X$ defines a map
$\Phi :\EuScript{V}_0 \times \RR \rightarrow \tV : (x,t) \mapsto \phi_t(x)$ 
which is a diffeomorphism on its image.
Let us introduce some notations: $\EuScript{V}_t = \Phi(\EuScript{V}_0 \times 
\{t\}) \subset \tV_{f(t)}$, its volume form will be denoted by $\omega_t$  
and  $\omega$  will be the
volume form of $\tV$. We know that $\omega_t = i_{X} \omega$. 
By definition of the divergence, we have 
$\d \omega_t = \d(i_{X} \omega)=-\div(X)\omega$, therefore
$\dd t \omega_t = -\div(X)\omega_t$. The computation
of $\div(X)$ gives
\begin{eqnarray*}
\div(X) & = & \div \left((2 \sqrt{\lambda_0} f)^{-1}\grad f  \right)\\
&=& (2 \sqrt{\lambda_0} f)^{-1} \div \grad f - i_{\grad \left(
(2 \sqrt{\lambda_0} f)^{-1} \right) } \grad f \\
&=& (2 \sqrt{\lambda_0} f)^{-1} \Delta f - (2 \sqrt{\lambda_0})^{-1}
\left\langle \grad \left(f^{-1}\right), \grad f \right\rangle \\
&=& -(2 \sqrt{\lambda_0})^{-1}
\left\langle -\frac{\grad f}{f^2}, \grad f \right\rangle = 
(2 \sqrt{\lambda_0}f^2)^{-1} |\grad f|^2 =2 \sqrt{\lambda_0}.
\end{eqnarray*}
We conclude that $\dd t \omega_t = -2 \sqrt{\lambda_0} \omega_t$, thus
$\omega_t = e^{-2 \sqrt{\lambda_0} t} \omega_0$
where $\omega_0$ is the $(n-1)$-volume form of $\tV_{f_0}$. As a result,
$\vol(\EuScript{V}_t) = e^{-2 \sqrt{\lambda_0} t} v_0$
where $v_0=\vol(\EuScript{V}_0)$ is a strictly positive constant.

As $(\Phi(\EuScript{V}_0\times \RR), \tg)$ is 
isometric to $(\EuScript{V}_0 \times \RR, g_t+\d t^2)$ 
where $(\EuScript{V}_0 \times \{t\}, g_t)$
is isometric to $(\EuScript{V}_t, \tg|_{\tV_{f(t)}})$, we can compute the following integral 
$$ \int_{\Phi(\EuScript{V}_0\times \RR)} f \omega = 
\int_\RR \int_{\EuScript{V}_t} f(t) \omega_t \d t =
 \int_\RR f_0 
e^{2 \sqrt{\lambda_0} t}e^{-2 \sqrt{\lambda_0} t} v_0 \d t
= \infty.$$
This is in contradiction with the fact that $f \in \L^1(\tV)$.
Hence 
we have proved the desired result  i.e.  $|\alpha|=0$.
\end{proof}

Finally we deduce a contradiction. Thanks to the equality 
(\ref{a}), Lemma \ref{AlphaZero} and ($a$), we write 
$|\tnabla \psi |=0=|\d |\psi||$, thus 
$|\psi|$ is constant. As $\psi$ belongs to $\L^2(\tV)$ with 
$\tV$ non-compact,
we deduct that $|\psi|=0$. This is impossible ($\psi$ is non-zero
by construction).
We conclude that if $\tV$ is non-compact, we cannot be in the equality case.

This ends the proof of Proposition \ref{FirstStep}.
\end{proof}

\begin{remark}
In the compact case we also have
$|\alpha|=0$. Therefore using the equality (\ref{a}),
 we have proved that in 
the equality case there exists 
a harmonic parallel 
non-zero spinor on $\tV$. This is a very strong constraint on the holonomy of
$\tV$.
The paper of M. Wang \cite{Wan89} is an illustration of this general 
philosophy.
\end{remark}

\subsection{Consequences and remarks}

Let $(V,g)$ satisfy the assumptions of Theorem 
\ref{AgenusCase}, then
using Theorem \ref{EqualityCase} we obtain the following
possiblities:
\begin{itemize}
\item[-] either $\scal(V,g) =0$ and then the equality holds in Theorem 
\ref{AgenusCase}, hence $\tV$ is compact. Moreover, as $\Ah_V[V] \neq 0$, we 
conclude that $\Ricci(V,g) \equiv 0$ (see \cite[VIII.6.]{Bou75} or 
\cite[Lem. $5.2.$]{KW75} for a proof).

\item[-] or $\inf \scal(V,g) < 0$ therefore we cannot be in the 
equality case, that is
$\inf \scal(V,g) < -4 \frac{n}{n-1} \lambda_0(\tV,\tilde{g})$.
\end{itemize}

We can compare this result with the Cheeger-Gromoll Theorem \cite{CG85}.
Our work proves again that a 
manifold $(V,g)$ with $\Ricci(V,g) \equiv 0$ and $\Ah_V[V] \neq 0$ has a finite
fundamental group (or equivalently, has a compact universal covering).

\section{Proof of Theorem \ref{EnlargeableCase}}
\label{Proof2}

Let $V^n$ be a closed Riemannian  manifold which
is enlargeable in the sense of Gromov-Lawson 
(see Definition \ref{Enlargeable}).

\subsection{Even dimensional case ($\dim V =n=2m$)}
\label{Proof2Even}

It is exactly the same idea as when M. Gromov and H.B. Lawson proved
the obstruction result cited in the introduction for enlargeable
manifolds of even dimension (\cite[$\S 5$]{GL83}, \cite{GL80}).
For sake of completeness we recall the method.

To prove Theorem \ref{EnlargeableCase} in the even dimensional case, 
we proceed \emph{by contradiction}.
Assume that the inequality (\ref{TheInequality}) does not hold i.e.
$$\inf \scal (V,g) > -4 \frac{n}{n-1} \lambda_0(\tV,\tilde{g}).$$
Then Theorem \ref{AgenusCase} asserts that $\Ah_V[V] = 0$.

In order to use Lemma \ref{UKato}, we must find ($\ZZ_2$-graded) 
Dirac operator on $V$ satisfying the condition of Propositition 
\ref{EvenCondition}.
As announced in Section \ref{Outline} we will use twisted Dirac operators
with unitary  bundles $Y$.
For such an operator, it is well known that $\ind({\Dirac^Y}^+) = 
\{\Ah(V) \ch(Y) \}[V]$. Hence we must construct bundles $Y$ such that 
$\{\Ah(V) \ch(Y) \}[V] \neq 0$.

\vskip1em
\paragraph{\textsf{Fundamental construction.}}
Let $f : V \rightarrow S^{2m}$ be a map of non-zero
degree. Consider the pullback bundle $Y$ via $f$ of a bundle 
$X$ on $S^{2m}$ (that we will specify afterwards).

\textit{Computation of the index of those twisted Dirac operators.}
In the case above we can compute the Chern character of the bundle $Y$.
We know that $\ch(Y) \in \H^{2*}(V)$ and $\ch(Y)= f^*\ch(X)$. Decomposing
$\ch(X) = \rank(X) + \ch_1(X)+ \cdots + \ch_m(X)$ where
$\ch_i(X) \in \H^{2i}(S^{2m})$, as we know the cohomology
of the sphere $S^{2m}$, we obtain
$$\ch(X) = \rank(X)+\ch_{m}(X) \textrm{ and therefore } 
\ch(Y) = \rank(X) + f^*\ch_{m}(X).$$  
Next considering the $\Ah$-class $\Ah(V)=1+\Ah_1(V)+ \cdots$ where
$\Ah_i(V) \in \H^{4i}(V)$, we check that
$$
\{\Ah(V) \ch(Y)\}[V] = \left\{
\begin{array}{l}
\{\rank(X) \Ah_{m/2}(V) + f^* \ch_m(X)\}[V] \textrm{ if $m$ is even,}\\
\{ f^* \ch_m(X)\}[V] \textrm{ if $m$ is odd.}\\
\end{array} \right.
$$
But $\{ \Ah_{m/2}(V)\}[V]=\Ah_V[V] = 0$, thus in the two cases
$$\{\Ah(V) \ch(Y)\}[V] =  \deg f \{\ch_m(X)\}[S^{2m}].$$

\textit{Choice of $X$ such that the index is non-zero.}
We have already chosen $f$ such that $\deg f \neq 0$. If we choose
$X$ equal to the positive half spinor bundle over $S^{2m}$, it is 
well known that $\{ch_m(X)\}[S^{2m}]\neq 0$. Therefore we have
constructed the desired bundle.

\vskip1em
\paragraph{\textsf{Use of enlargeability.}}
Note that, as $Y=f^* X$, we get $\| R^Y\| \leq \|\d f\|^2 \|R^X\|$.
Using enlargeability, for all $\epsilon >0$, we can find
a finite Riemannian spin covering $\tV_\epsilon$ of $V$ and a map
$f_\epsilon :\tV_\epsilon \rightarrow S^{2m}$ of
non-zero degree such that the bundle $Y_\epsilon$ constructed
as above satisfies
$$\|R^{Y_\epsilon}\| \leq \epsilon \;\;\textrm{ and }\;\;
\{ \Ah(\tV_\epsilon) \ch(Y_\epsilon) \} [\tV_\epsilon] \neq 0.$$ 
(Since the $\Ah$-genus is multiplicative under coverings, 
the $\Ah$-genus of $V$ and $\tV_\epsilon$ are both zero.)
In fact we consider the following diagram:

$$\xymatrix@1{
Y_\epsilon = f_\epsilon^*(X)  \ar[d] & X \ar[d]\\
\tV_\epsilon \ar[r]^-{f_\epsilon} \ar@{.>}[d] & S^{2m} \\
V 
}$$

Then we construct the twisted Dirac operator $\Dirac_\epsilon$ on
$\tV_\epsilon$ acting on sections of the bundles $\Spin \otimes Y_\epsilon$.
Let $\tDirac_\epsilon$ be its lift on $\tV$. We have $\ind_\G 
(\tDirac_\epsilon^+) = \ind(\Dirac^+_\epsilon) = \{ \Ah(\tV_\epsilon)
 \ch(Y_\epsilon) \}[\tV_\epsilon]\neq 0$. Using Proposition
\ref{EvenCondition}, we deduce that $0$ belongs to the point spectrum of
$\tDirac_\epsilon$. Thus we can apply Lemma \ref{UKato} to the manifold
$\tV_\epsilon$ to obtain
$$\inf \scal(V,g) = \inf \scal(\tV_\epsilon,\tilde{g}_\epsilon)
\leq -4 \frac{n}{n-1} \lambda_0(\tV, \tilde{g})+ \alpha_n \epsilon.$$
Letting $\epsilon \rightarrow 0$, we get a contradiction. Consequently the
inequality (\ref{TheInequality}) holds.

\subsection{Odd dimensional case ($\dim V = n=2m-1$)}
\label{Proof2Odd}

Since Section \ref{Index} we know that if we use a family of 
Dirac operators in order to find an operator $\tD$ on $\tV$ such
that $0$ belongs to the spectrum of $\tD$, we must satisfy  
Theorem \ref{OddCondition}: the spectral flow must be non-zero. 
In particular, we will look for 
a family of twisted Dirac operators $\Dirac^{Y_u}$ such 
that its spectral flow is  non-zero.
In fact it is exactly the same idea as in the even dimensional case:
we pull back a bundle from $S^{2m-1} \times [0,1]$ to $V \times [0,1]$. 
Then, the key
Lemma \ref{UKato} permits to  conclude using enlargeability. 

\vskip1em
\paragraph{\textsf{Fundamental construction.}}
\label{Construction}

We borrow this construction from M. Atiyah in \cite[p. 258-259]{Ati85}.
Let $h : S^{2m-1} \rightarrow \Gl_N(\CC)$ 
(that we will choose later)
where we take $N$ sufficiently large.
Fix a trivialization on the trivial bundle $X_0=S^{2m-1} \times \CC^N$ and
consider $h$ as a 
multiplication operator on this bundle. 
We denote $\nabla_0$ the trivial connection on $X_0$ with respect to the fixed 
trivialization and 
$\nabla_1=h^*\nabla_0$. Consider now the linear family of connections
$\nabla_u = u \nabla_1 + (1-u) \nabla_0$ joining $\nabla_0$ to its gauge
transform $\nabla_1$. We also deduce a (trivial) bundle $X$ 
on $S^{2m-1} \times [0,1]$
with connection $\nabla= \nabla_u + \pdd{u}$ and a fixed trivialization.

With a map $f : V \rightarrow S^{2m-1}$ of non-zero degree 
 we construct the 
corresponding family of bundles $(Y_0, A_u) = (f^*X_0, f^* \nabla_u)$ on $V$,
the (trivial) bundle $Y$ on $V \times [0,1]$ with connection
$A = A_u+ \pdd{u}$ and the corresponding family of 
twisted Dirac
operator $\D_u = \Dirac^{Y_u}$.
We want to find a family $Y_u$ (i.e. maps $f$ and $h$) such that
the spectral flow is non-zero.

\textit{Computation of the spectral flow of this family.}
In the case of twisted Dirac operators $\Dirac^{Y_u}$ and  
of the operator $\mathcal{D}=\pdd u -\Dirac^{Y_u}$
(take care on the convention on $\mathcal{D}$, see Remark \ref{Convention}),
we know that
$A(1) = - \{ \Ah(V \times [0,1]) 
\ch(Y,A)\} [V \times [0,1]]$. Moreover, thanks
Lemma \ref{sflow}, the spectral flow is equal to  
$-A(1)$.

Let notice that, on a manifold \emph{with boundary},
the Chern character depends 
on the connection. Thus even if a bundle is trivial,
it can have a non-trivial Chern character.

By construction $\ch (Y,A) = (f \times \id)^* \ch(X,\nabla)$. 
Moreover, the fact that $ \ch (X,\nabla) \in \H^{2*}(S^{2m-1} \times [0,1])$
and the K\"unneth formula assert 
$$\ch(X,\nabla) = \rank(X) + \ch_{m}(X,\nabla), \textrm{ hence }
 \ch(Y,A) = \rank(X) + (f \times \id)^* \ch_{m}(X,\nabla).$$ 
We have also that $\Ah(V \times [0,1]) = \Ah(V) \in \H^{4*}
(V \times [0,1])$. Thus
{\setlength\arraycolsep{2pt}
\begin{eqnarray*}
\sflow(1) & = & \{ \Ah(V \times [0,1]) \ch(Y,A)\}[V \times [0,1]]\\
& = & \{ \Ah(V) \left(\rank (X) + (f\times \id)^* \ch_m(X,\nabla)\right)\}
[V \times [0,1]]\\
& = &   \rank(X) \Ah_V[V] + \deg f \{ch_m(X,\nabla)\}
[S^{2m-1} \times [0,1]].
\end{eqnarray*}}
For an odd dimensional manifold, we always have a zero $\Ah$-genus.
It remains to know how $\{ch_m(X,\nabla)\}
[S^{2m-1} \times [0,1]]$ depends on $h$.

To begin, in the above fixed trivialization, 
the connection $1$-form of $(X_0, \nabla_1=h^*\nabla_0)$ is given by
$\omega^1 = h^{-1} \d h$. 
Thus the connection $1$-form of  the bundle 
$(X_0,\nabla_u=u\nabla_1+(1-u)\nabla_0)$ is $\omega^u = u h^{-1} \d h + 
(1-u)\omega_0$ but $\nabla_0$ is trivial in the fixed trivialization
hence $\omega_0=0$.
Therefore the connection $1$-form of  the bundle 
$(X,\nabla)$ at the point $(s,u) \in S^{2m-1}\times [0,1]$ is 
$$\omega_{(s,u)}= \omega^u_s= u  h^{-1} \d h .$$
Now we calculate the connection $2$-form of  the bundle 
$(X,\nabla)$ remarking that $d(h^{-1})=-h^{-1}\d h h^{-1}$ 
$$\Omega= \d u  \wedge  h^{-1} \d h - (u+u^2) h^{-1} \d h \wedge h^{-1} \d h.$$
Recall that a differential form representing  $\ch_{m}(X)$ is  given by
$$\left( \frac{i}{2 \pi} \right)^m \frac{1}{m !} \Tr(\Omega^m).$$
If $h^{-1} \d h$ appears an even number of times, 
the trace is zero. If $\d u$ appears
twice, the form is zero. Therefore 
{\setlength\arraycolsep{2pt}
\begin{eqnarray*}
 \left( \frac{i}{2 \pi} \right)^m \frac{1}{m !} \Tr(\Omega^m)
 & = & C \Tr( (u+u^2)^{m-1}
\d u \wedge \underbrace{h^{-1} \d h \wedge \cdots \wedge h^{-1} \d h}_{2m-1})\\
& = & C (u+u^2)^{m-1}\d u \wedge \Tr(h^{-1} \d h \wedge 
\cdots \wedge h^{-1} \d h)
\end{eqnarray*}}
consequently
$$\{ch_m(X,\nabla)\}[S^{2m-1} \times [0,1]] = C
\int_0^1 \! \! \! \! (u+u^2)^{m-1}\d u
\int_{S^{2m-1}} \! \! \! \! \! \Tr(h^{-1} \d h 
\wedge \cdots \wedge h^{-1} \d h)$$
where $C$ is a non-zero constant.

\textit{Choice of $h$ such that the spectral flow is non-zero.}
We have already chosen $f$ such that $\deg f \neq 0$ and
we have trivially that $\int_0^1(u+u^2)^{m-1}\d u \neq 0$. Thus we must
choose $h$  
such  that  
$\int_{s^{2m-1}}\Tr(h^{-1} \d h \wedge \cdots \wedge h^{-1} 
\d h) \neq 0$.

For constructing $h$, 
we begin with choosing $J_0, \cdots, J_{2m}$ in 
$\operatorname{Aut}(V)=\Gl_{2N}(\CC)$
 which satisfy the identities 
$$J_\alpha J_\beta = -J_\beta J\alpha \;\textrm{ if }\; \alpha \neq \beta \;\;
\textrm{ and } \;\;(J_\alpha)^2 = \Id.$$
($N$ must be sufficiently large so that $\CC^N$ is such a Clifford module,
for example $N=2^m$).
Let $J$ be $i^m J_1J_2 \cdots J_{2m}$. As $J^2=\Id$, if we denote
 $V_0 = \Ker(J-\Id)$ and $V_1= \Ker(J+Id)$, we have $V=V_0 \oplus V_1$, each
automorphism $J_\alpha$ maps $V_0$ onto $V_1$ 
(note that $J_\alpha J=-JJ_\alpha$) and
$\dim V_0 = \dim V_1 = N$. We also remark that for an even number $k <m$ 
of distinct $\alpha_i$ we have $\Tr((J_{\alpha_1}\cdots
J_{\alpha_{2k}})_{|V_0})=0$. 
Now we define $h : S^{2m-1} \rightarrow \operatorname{Aut}(V_0)=\Gl_{N}(\CC)$
by $h(x_1,\cdots,x_{2m})=  J_1(x_1J_1+x_2J_2+ \cdots +x_{2m}J_{2m})$ 
restricted to 
$V_0$. By noticing that $\d (h^{-1})=- h^{-1} \d h h^{-1}$, that
$(x_1J_1+\cdots+x_{2m}J_{2m})^2=  \Id$ and
that $(J_\alpha \d x_\alpha) \wedge (J_\beta \d x_\beta)= 
(J_\beta \d x_\beta) \wedge
(J_\alpha \d x_\alpha)$ we can prove that 
$$\Tr(\underbrace{h^{-1} \d h \wedge \cdots \wedge 
h^{-1} \d h}_{2m-1}) = (2m-1)! N 
\sum_{j=1}^{2m} (-1)^j x_j \d x_1 \wedge \cdots \wedge \widehat{\d x_j} \wedge \cdots \wedge
\d x_{2m}.$$
i.e. $\Tr(h^{-1} \d h \wedge \cdots \wedge h^{-1} \d h)$ 
is proportional to the volume form
of the sphere $S^{2m-1}$.
We conclude that  $\int_{S^{2m-1}}\Tr(h^{-1} \d h \wedge \cdots 
\wedge h^{-1} \d h) \neq 0$.
Moreover $h$ define a non-zero element in $\pi_{2m-1}(\Gl_{N}(\CC))= \ZZ$.

\vskip1em
\paragraph{\textsf{Use of enlargeability.}}
\label{Conclusion}

It is the same as in the even dimensional case. For a family constructed 
as above, using the key Lemma \ref{UKato}, we have
$$\inf \scal(V,g) \leq -4 \frac{n}{n-1} \lambda_0(\tV,\tilde{g}) +
\alpha_n \|R^{Y_{u_0}}\| \;\;\textrm{ for some } u_0 \in [0,1].$$
Recall that $Y_u = f^*(X_u)$ where $f : V \rightarrow S^{2m-1}$
is of non-zero degree.
Thus $\|R^{Y_u}\| \leq \|\d f\|^2\|R^{X_u}\|$. But by construction, we 
know that $\|R^{X_u}\| \leq C$ 
where $C$ is a constant independent of $u$.
Therefore we have
 $$\|R^{Y_u}\| \leq C \|\d f\|^2.$$
If we replace $V$ by a Riemannian finite covering (using enlargeability),
 this does not change the result and we can choose $f$ such that $\|R^{Y_u}\|$ 
is as small as we want, for example $\|R^{Y_u}\| \leq \epsilon$ with 
$\epsilon >0$. Letting $\epsilon \rightarrow 0$ in Lemma \ref{UKato}, 
we  conclude as in the even dimensional case.

\section{Generalizations and Remarks} 
\label{Generalization}

\subsection{Unification of Theorem \ref{AgenusCase} and Theorem 
\ref{EnlargeableCase}}

We can use a more general notion of enlargeability
(see \cite[$\S 5 \& 6$]{GL83}), more precisely:

\begin{Edefinition}
A closed Riemannian $n$-manifold $(V,g)$ is said to be 
\emph{area-enlargeable in dimension $p$} if
given any $\epsilon >0$, there exists a finite Riemannian covering 
$\tV_\epsilon$ of $V$ which is spin and a $C^1$-map $f_\epsilon :
\tV_\epsilon \rightarrow S^p$ which is $\epsilon$-contracting on $2$-vectors
(i.e. if $\|f_* \phi \| \leq \epsilon \|\phi\|$ for all $2$-vectors 
$\phi$ on $\tV_\epsilon$) and of non-zero $\Ah$-degree (i.e. the $\Ah$-genus
$\Ah(f^{-1}(y))$ of the inverse image of each regular value $y$ is non-zero).
\end{Edefinition}

We remark that for a closed Riemannian $n$-manifold $(V,g)$ which is
area-enlargeable in dimension $p$, we have necessary $p \equiv n \mod 4$.
The case $p=0$ correspond to the non-zero $\Ah$-genus case. The case
$p=n$ correspond to the enlargeable manifold as defined in the introduction.
We use \emph{area-contracting map} because we notice that in the proof
of Theorem \ref{EnlargeableCase} we only need the contraction of $2$-vectors
(and not necessary of vectors themselves). Thus we can give the following
unified statement:

\begin{theo}
\label{Unified}
Let $(V,g)$ be a closed Riemannian manifold of dimension $n$. 
If $V$ is area-enlargeable in dimension $p$, 
then the  inequality (\ref{TheInequality}) holds.
\end{theo}

\subsection{On $\Karea$}

In the paper \cite{Gro96}, M. Gromov sets the notion of enlargeability
in a broader perspective. He defines  a Riemannian invariant, the $\Karea$. 
Roughly speaking, $(V,g)$ has a big $\Karea$ if it carries a (non trivial !)
vector bundle with small curvature.
More precisely:

\begin{Edefinition}
Let $(V^{2m},g)$ be a closed even-dimensional oriented Riemannian manifold.

A bundle $X$ over $V$ is said \emph{homologically significant} if at least
one characteristic (Chern) number of $X$ does not vanish.

Let 
$$C(V,g)= \inf \|R^X\|$$
be the infimum taken on all hermitian homologically significant
bundles $X$. We define the $\Karea$ to be the inverse of this number:
$$\Karea(V,g)=C(V,g)^{-1}.$$
\end{Edefinition}

It follows immediately that if $V$ is enlargeable, $\Karea(M,g)$ is 
infinite for any $g$. In fact the property for the $\Karea$ to be infinite is,
like enlargeability, a topological property.
For details on this new Riemannian invariant we refer to the 
original paper \cite[$\S 4 \& 5$]{Gro96}.

Then M. Gromov proves a theorem which can replace the Gromov-Lawson Theorem 
 on enlargeable manifolds (\cite{GL80}). 
It states that for  a closed Riemannian spin manifold $V$, 
if
$\inf \scal(V,g) \geq \epsilon^{2}$ then
$\Karea(V,g) \leq C \epsilon^{-2}$  if $V$ is even dimensional
and  $\Karea(V\times S^1,g\times \mathrm{can}) \leq C \epsilon^{-2}$ 
if $V$ is odd dimensional ;
where $C$ is a constant which depends only on the dimension of $V$.
As a consequence: a closed even-dimensional spin 
Riemannian  manifold of infinite $\Karea$ does not carry any  
metric of strictly positive scalar curvature. 
Therefore
the proof of Theorem \ref{EnlargeableCase} 
can be adapted to obtain
the following:

\begin{theo} 
\label{KareaCase}
Let $(V,g)$ be a closed Riemannian spin manifold of dimension $n$. 
If $V$ is even dimensional and $\Karea(V) = \infty$ or if $V$ is odd
dimensional and $\Karea(V \times S^1) = \infty$, 
then the  inequality (\ref{TheInequality}) holds.
\end{theo}

\begin{proof}
To prove this theorem,
it suffices to remark that the assumptions on $\Karea$ allow
to construct the desired bundles in the even dimensional as
in the odd dimensional cases and then to use Lemma \ref{UKato}.

If the $\Karea$ of
a closed even-dimensional manifold $V$ is infinity, for every $\epsilon >0$
there exists a bundle $Y_\epsilon$ on $V$ such that
$\|R^{Y_\epsilon}\| \leq \epsilon$ and $\{ \widehat{A}(V) \ch(Y_\epsilon)\}
\neq 0$. Afterwards, the proof is exactly the same as for an enlargeable 
even-dimensional manifold.
 
In the odd dimensional case, recall that  the spectral flow
only depends on Atiyah-Singer index Theorem on the compact manifold
$V \times S^1$ (see Remark \ref{Bismut-Freed}). Hence it suffices 
to replace $V$ by $V \times S^1$ in the previous statement to obtain the 
result.
\end{proof}

\subsection{Consequences}

In \cite[$\S 5 \frac{1}{2}$]{Gro96}, M.Gromov gives
the following corollary which is a very interesting
non-approximation interpretation of the 
inequality of K. Ono (or {\it a fortiori} ours) 
and we give here a concise proof.

\begin{corollary}[M. Gromov] 
\label{Gromov}
Let $(V,g_0)$ be a closed Riemannian manifold which
is area-enlargeable in dimension $p$ and  with a non amenable 
fundamental group.
Then, there exist a constant 
$\sigma_0 = \sigma_0(g_0) > 0$ such that $g_0$ admits no $C^0$-approximation
by $C^2$-metrics $g$ on $V$ with $\scal(V,g) \geq -\sigma_0$.
\end{corollary}

\begin{proof}
The Riemannian manifold $(V,g)$ satisfies the assumptions of Theorem
\ref{Unified}, consequently 
$$\inf \scal(V,g) \leq -4 \frac{n}{n-1} \lambda_0(\tV, \tilde{g}).$$
A basis of neighbourhoods of $g_0$ is given by the sets
$\mathcal{V}_{ab}=\{g ,\;\; a^2 g_0 \leq g \leq b^2 g_0\}$ where $0<a^2<1<b^2$.
Then for all $g \in \mathcal{V}_{ab}$ and for all $k \in \NN$ we have 
$\frac{a^n}{b^{n+2}} \lambda_k(\tilde{g}_0) 
\leq \lambda_k(\tilde{g})
\leq \frac{b^n}{a^{n+2}} \lambda_k(\tilde{g}_0)$ from the Courant-Fisher
min-max characterization of the eigenvalues of the Laplace operator.
This implies that the function $g \mapsto \lambda_0(\tV, \tilde{g})$ is 
continuous therefore
there exists a $C^0$-neighbourhood $\mathcal{V}$ of $g_0$ such that for all
metric $g$ in $\mathcal{V}$
$$|\lambda_0(\tV, \tilde{g})-\lambda_0(\tV, \tilde{g_0})| <\epsilon.$$

As the fundamental group is non amenable, we know that
$\lambda_0(\tV, \tilde{g_0})\neq 0$
(\cite{Bro81}).  
Letting $\sigma_0(g_0)= 4 \frac{n}{n-1}( \lambda_0(\tV, \tilde{g_0}) 
-\epsilon)>0$ (for $\epsilon$ sufficiently small) then
for all
metric $g$ in $\mathcal{V}$, 
$$\inf \scal(V,g) \leq -4 \frac{n}{n-1} 
\lambda_0(\tV, \tilde{g}) <  -4 \frac{n}{n-1} 
(\lambda_0(\tV, \tilde{g_0})-\epsilon) = -\sigma_0(g_0).$$
Consequently we cannot approximate $g_0$ by a $C^0$-sequence of $C^2$-metrics 
such that
 $\inf \scal(V,g) \geq  -\sigma_0(g_0)$.
\end{proof}
\section*{Acknowledgements}
I would like to thank Jacques Lafontaine 
for his advice and encouragement. 
I am also grateful to Marc Herzlich for explaining the refined Kato inequality
to me and  for his
many useful comments and careful reading of the manuscript.
Finally, I would like to thank Xiaonan Ma for pointing out  
the paper of Bismut-Freed, as well as Gilles Carron for helpful
communications. I also wish to thank M. Min-Oo for 
his series of lectures on
$\Karea$ at the CRM (Centre de Recerca Matem\`atica) during Summer
$2001$.
\end{document}